\documentclass[11pt]{amsart}
\usepackage{geometry}                
\usepackage[parfill]{parskip}    
\usepackage{graphicx}
\usepackage{amssymb}
\usepackage{epstopdf}
\newtheorem{theorem}{Theorem}[section]
\newtheorem{corollary}{Corollary}

\newtheorem{lemma}[theorem]{Lemma}

\theoremstyle{definition}
\newtheorem{definition}[theorem]{Definition}

\newcommand{\R}{\mathbb{R}}
\newcommand{\N}{\mathbb{N}}

\newcommand{\Z}{\mathbb{Z}}

\title[Carrying simplex]
      {On existence and uniqueness of a modified carrying simplex for discrete Kolmogorov systems}

\author[Zhanyuan Hou]{}

\subjclass{Primary: 37B25; Secondary: 37C70, 34D23, 34D05.}
 \keywords{discrete competitive models, retrotone maps, carrying simplex, existence and uniqueness, dominant species, vanishing species}

 \email{z.hou@londonmet.ac.uk}

\begin{document}

\maketitle

\centerline{\scshape Zhanyuan Hou }
\medskip
{\footnotesize
 \centerline{School of Computing and Digital Media, London Metropolitan University,}
   \centerline{166-220 Holloway Road, London N7~8DB, UK}
} 

\medskip


\begin{abstract}
For a $C^1$ map $T$ from $C =[0, +\infty)^N$ to $C$ of the form $T_i(x) = x_if_i(x)$, the dynamical system $x(n) =T^n(x)$ as a population model is competitive if $\frac{\partial f_i}{\partial x_j}\leq 0$ $(i\not= j)$. A well know theorem for competitive systems, presented by Hirsch (J. Bio. Dyn. 2 (2008) 169--179) and proved by Ruiz-Herrera (J. Differ. Equ. Appl. 19 (2013) 96--113) with various versions by others, states that, under certain conditions, the system has a compact invariant surface $\Sigma\subset C$ that is homeomorphic to $\Delta^{N-1} =\{x\in C: x_1+ \cdots + x_N=1\}$, attracting all the points of $C\setminus\{0\}$, and called carrying simplex. The theorem has been well accepted with a large number of citations. In this paper, we point out that one of its conditions requiring all the $N^2$ entries of the Jacobian matrix $Df = (\frac{\partial f_i}{\partial x_j})$ to be negative is unnecessarily strong and too restrictive. We prove the existence and uniqueness of a modified carrying simplex by reducing that condition to requiring every entry of $Df$ to be nonpositive and each $f_i$ is strictly decreasing in $x_i$. As an example of applications of the main result, sufficient conditions are provided for vanishing species and dominance of one species over others.
\end{abstract}

\textbf{Note.} This paper has been accepted for publication in Journal of Difference Equations and Applications.


\section{Introduction}\label{Sec1}
In this paper, we are concerned with the global asymptotic behaviour of the discrete dynamical system
\begin{equation}\label{e1}
x(n) = T^n(x), \quad x\in C,\quad n\in \N,
\end{equation}
where $C = \R^N_+ = [0, +\infty)^N$, $\N = \{0, 1, 2, \ldots\}$ and the map $T: C \to C$ has the form
\begin{equation}\label{e2}
T_i(x) = x_if_i(x), \quad i\in I_N = \{1, 2, \ldots, N\}
\end{equation}
and $f\in C^1(C, C)$ with $f_i(x)>0$ for all $x\in C$ and $i\in I_N$. System (\ref{e1}) is a typical mathematical model for the population dynamics of a community of $N$ species, where each $x_i(n)$ represents the population size or density at time $n$ (at the end of $n$th time period), and the function $f_i(x)$ denotes the per capita growth rate, of the $i$th species. If $\frac{\partial f_i}{\partial x_j} \leq 0$ for all $i, j\in I_N$ with $i\not= j$, then increase of the $j$th population reduces the per capita growth rate of the $i$th species, so (\ref{e1}) models the population dynamics of a community of competitive species. 

System (\ref{e1}) and its various particular instances as models have attracted huge interests from researchers in the last two decades. One of the important and influential developments is the existence of a {\it carrying simplex} $\Sigma \subset C$: a compact invariant hypersurface homeomorphic to $\Delta^{N-1} =\{x\in C: x_1 + \cdots +x_N =1\}$ such that every trajectory except the origin is asymptotic with a trajectory in $\Sigma$. Since $\Sigma$ attracts all the points of $C\setminus\{0\}$, the dynamics of (\ref{e1}) on $C$ is essentially described by the dynamics on $\Sigma$. The carrying simplex theory was originally established by Hirsch \cite{Hir1} (see \cite{Hou1} for latest update) for competitive Kolmogorov systems of differential equations. Since then the idea of a carrying simplex for discrete systems gradually appeared in literature (see \cite{WaJi1}, \cite{WaJi2}, \cite {JiMiWa} for example). But a more accepted theorem for existence and uniqueness of a carrying simplex for (\ref{e1}) was given by Hirsch \cite{Hir2} without proof. Then Ruiz-Herrera \cite{Her} presented a more general theorem covering Hirsch's result with a complete proof.

For any $x, y\in C$, we write $x\leq y$ or $y \geq x$ if $x_i\leq y_i$ for all $i\in I_N$; $x<y$ or $y>x$ if $x\leq y$ but $x\not= y$; $x\ll y$ or $y\gg x$ if $x_i < y_i$ for all $i\in I_N$. The map $T$ given by (\ref{e2}) is said to be {\it retrotone} in a subset $X\subset C$ if for any $x, y\in X$, $T(x) < T(y)$ implies $x_i < y_i$ for all $i\in I(y) = \{j\in I_N: y_j\not= 0\}$. Let $[0, r] = \{x\in C: 0\leq x\leq r\}$. The theorem below is Theorem 6.1 in \cite{Her}.

\begin{theorem}\label{The1.1}
	Assume that $T$ with $T([0, r]) \subset [0, r]$ for some $r\gg 0$ satisfies the following conditions:
	\begin{itemize}
		\item[(i)] For each $i\in I_N$, the map $T$ restricted to the positive half $x_i$-axis has a fixed point $q_ie_i$ with $q_i>0$, $e_i$ the $i$th standard unit vector and $q\ll r$.
		\item[(ii)] $T$ is retrotone and locally one to one in $[0, r]$.
		\item[(iii)] For any $x, y\in [0, r]$, if $T(x) < T(y)$ then, for each $j\in I_N$, either $x_j=0$ or $f_j(x) >f_j(y)$.
	\end{itemize}
	Then the map admits a carrying simplex $\Sigma$.
\end{theorem}

Note that Theorem \ref{The1.1} can be only applied to the system restricted to the space $[0, r]\subset C$ if no condition for $T$ on $C\setminus[0, r]$ is provided. However, if for any compact set $S\subset C$ there is a $k\in\N$ such that $T^k(S) \subset [0, r]$, then Theorem \ref{The1.1} can be applied directly to the system on $C$.

When $f$ on $C$ is a $C^1$ map, $T$ is also a $C^1$ map with Jacobian matrix
\begin{equation}\label{e3}
DT(x) = \textup{diag}(f_1(x), \ldots, f_N(x))(I - M(x)),
\end{equation}
where $I$ is the identity matrix and 
\begin{equation}\label{e4}
M(x)= (M_{ij}(x))= \left(-\frac{x_i}{f_i(x)}\frac{\partial f_i}{\partial x_j}(x)\right)_{N\times N}.
\end{equation}
Then, by Lemma 4.1, Corollary 6.1 and Remark 6.4 in \cite{Her}, Theorem \ref{The1.1} has the following version with easily checkable conditions.

\begin{theorem}\label{The1.2}
	Assume that $T$ satisfies the following conditions:
	\begin{itemize}
		\item[(i)] For each $i\in I_N$, the map $T$ restricted to the positive half $x_i$-axis has a fixed point $q_ie_i$ with $q_i>0$, $e_i$ the $i$th standard unit vector and $q\ll r$ for some $r\in C$.
		\item[(ii)] All entries of the Jacobian $Df$ are negative.
		\item[(iii)] The spectral radius of $M(x)$ satisfies $\rho(M(x)) <1$ for all $x\in [0, q]\setminus\{0\}$.
	\end{itemize}
	Then the map admits a carrying simplex $\Sigma$.
\end{theorem}

A more user-friendly variation of Theorem \ref{The1.2} given by Jiang and Niu \cite[Theorem 3.1]{JiNi2} is the above theorem with simply a replacement of condition (iii) by (iii$)'$ below:
\begin{itemize}
	\item[(iii$)'$] For each $x\in [0, q]\setminus\{0\}$ with $I(x) = \{j\in I_N: x_j>0\}$, either 
	\[
	f_i(x) + \sum_{j\in I(x)}x_j\frac{\partial f_i}{\partial x_j}(x)>0 \quad \forall i\in I(x)
	\] 
	or 
	\[
	f_i(x) + \sum_{j\in I(x)}x_i\frac{\partial f_i}{\partial x_j}(x)>0 \quad \forall i\in I(x).
	\]
\end{itemize}

A carrying simplex $\Sigma$ has the important and interesting features: compact, invariant, unordered ($p\leq q$ implies $p=q$ for $p, q\in \Sigma$), homeomorphic to $\Delta^{N-1}$ by radial projection, and attracting all the points of  $C\setminus\{0\}$. Therefore, if (\ref{e1}) admits a carrying simplex, the dynamics of the system on the $N$-dimensional space $C$ is essentially described by the dynamics on this $(N-1)$-dimensional hypersurface $\Sigma$. Due to these features of $\Sigma$, Theorem \ref{The1.2} is phenomenal and lays the foundation for further investigations. There are a large number of applications of this theorem, the following are just a few examples.

Ruiz-Herrera \cite{Her} investigated exclusion and dominance utilizing the existence of a carrying simplex. Jiang and Niu \cite{JiNi1, JiNi2} and Gyllenberg et al. \cite{GyJiNiYa1, GyJiNiYa2} dealt with some well known three-dimensional competitive models. Based on the existence of a carrying simplex, they classified the systems into 33 topologically equivalent classes and gave a phase portrait on $\Sigma$ for each class. Jiang, Niu and Wang \cite{JiNiWa} studied heteroclinic cycles via carrying simplex. Balreira et al. \cite{BaElLu} and Gyllenberg et al. \cite{GyJiNi3} provided criteria for global stability of an interior fixed point based on the existence of a carrying simplex. Baigent \cite{Bai1, Bai2} investigated the geometric feature of a carrying simplex and found conditions for $\Sigma$ to be convex. Baigent and Hou \cite{BaHo} and Hou \cite{Hou} provided split Lyapunov function method and geometric method for global stability. Although these methods were not based on the existence of a carrying simplex, comments and comparisons with those using carrying simplex were made there.

We note that condition (ii) in Theorem \ref{The1.2} is very restrictive; it excludes the possibility of applying the theorem to systems with some zero entries of $Df$. But actually, condition (ii) is too strong and unnecessary, a compact invariant set attracting all the points of $C\setminus\{0\}$ with most of the features of a carrying simplex may still exist even if $\frac{\partial f_i}{\partial x_j} =0$ for some distinct $i, j\in I_N$.

The aim of this paper is to prove the existence and uniqueness of a modified carrying simplex under a much weaker condition than (ii): instead of (ii) requiring all $N^2$ entries of $Df$ to be negative, we require each entry of $Df$ to be nonpositive, with each $f_i$ strictly decreasing in $x_i$, on a compact set. We shall present the main results in section 2 and leave the proofs to section 5. In section 3, we present some results on dominant species and vanishing species as an application of the main results. In section 4, we deal with some known models as examples. We finally conclude the paper in section 6.

\section{Notation and main results}\label{Sec2}
For $C = \R^N_+$ we let $\dot{C} = \{x\in C: \forall i\in I_N, x_i >0\}$ and $\partial C = C\setminus \dot{C}$. Then $\dot{C}$ is the interior of $C$ and $\partial C$ is the boundary of $C$. The part of $\partial C$ restricted to the $i$th coordinate plane and the part restricted to the positive half $x_i$-axis are denoted by $\pi_i$ and $X_i$ respectively, i.e.
\begin{eqnarray*}
	\pi_i &=& \{x\in C: x_i =0\}, i\in I_N,\\ 
	X_i &=& \{x\in C: x_i>0, \forall j\in I_N\setminus\{i\}, x_j =0\}, i\in I_N.
\end{eqnarray*}
Denote the $i$th standard unit vector by $e_i$, i.e. the $i$th component of $e_i$ is 1 and others are 0. For any nonempty subset $I\subset I_N$, define
\begin{eqnarray*}
	C_I &=& \{x\in C: \forall j\in I_N\setminus I, x_j =0\}, \\
	\dot{C}_I &=& \{x\in C_I: \forall i\in I, x_i >0\}.
\end{eqnarray*}
For any $x, y\in C_I$, by writing $x\leq_I y$ we mean $x_i \leq y_i$ for all $i\in I$; we write $x<_I y$ if $x\leq_I y$ but $x\not= y$; and we write $x\ll_I y$ if $y-x \in \dot{C}_I$. We may also use $y\geq_I x$, $y>_I x$ and $y\gg_I x$ for $x\leq_I y$, $x<_I y$ and $x\ll_I y$ respectively. If $I = I_N$, we simply drop the subscript ``$I$'' from these inequalities. For any $a, b\in C$ with $a\leq b$, we let $[a, b] = \{x\in C: a\leq x\leq b\}$. Then $[a, b]$ is a $k$-dimensional cell if $b-a$ has exactly $k$ positive components. For each $x\in C$, the positive limit set $\omega(x)$ of $T^n(x)$ is defined by
\[
\omega(x) = \bigcap^{\infty}_{n=1}\overline{\{T^k(x): k\geq n\}},
\]
where $\overline{A}$ denotes the closure of any set $A$. If $T$ is invertible and $T^{-n}(x)$ exist for all $n\in\N$, the negative limit set $\alpha(x)$ is defined by
\[
\alpha(x) = \bigcap^{\infty}_{n=1}\overline{\{T^{-k}(x): k\geq n\}}.
\] 
Also, the whole trajectory of $x$ is denoted by $\gamma(x) = \{T^n(x): n\in \Z\}$. 

Suppose a simply connected closed set $S\subset C\setminus\{0\}$ is an $(N-1)$-dimensional hypersurface which divides $C$ into three mutually exclusive subsets $S^-$, $S$ and $S^+$ with $0\in S^-$ and $C = S^-\cup S\cup S^+$. A point $p\in C$ is said to be {\it below} ({\it on} or {\it above}) $S$ if $p\in S^-$ ($S$ or $S^+$). For any nonempty subset $S_0\subset C$, $S_0$ is said to be {\it below} ({\it above}) $S$ if $S_0\subset S^-\cup S$ ($S\cup S^+$); $S_0$ is said to be {\it strictly below} ({\it strictly above}) $S$ if $S_0\subset S^-$ ($S^+$). 

Let $B$ be either $C$ or a positively invariant $[0, r]$ for some $r\in\dot{C}$. For convenience, we define the concept of a modified carrying simplex as follows.

\begin{definition}
	A nonempty set $\Sigma\subset B\setminus\{0\}$ is called a modified carrying simplex of (\ref{e1}) if $\Sigma$ meets the following requirements.
	\begin{itemize}
		\item[(i)] $\Sigma$ is compact, invariant and homeomorphic to $\Delta^{N-1}$ by radial projection.
		\item[(ii)] $\Sigma$ attracts all the points of $B\setminus\{0\}$, i.e. $\omega(x)\subset \Sigma$ for each $x\in B\setminus\{0\}$.  
	\end{itemize}
	Moreover, if $x$ is below $\Sigma$ with a nonempty support $I(x)\subset I_N$, then there is a $y\in\Sigma$ with $I(y)=I(x)$ such that $\lim_{n\to +\infty}(T^n(x)-T^n(y)) =0$.
\end{definition}

Note that the ``unordered'' property of $\Sigma$ is not mentioned in the above definition. We shall see in Remark 2.1 (f) below that the unordered property of $\Sigma$ here is slightly different from that for carrying simplex in Hirsch \cite{Hir2}, Ruiz-Herrera \cite{Her} and the literature. But the main difference between modified carrying simplex and the carrying simplex in literature is that the latter requires every trajectory in $B\setminus\{0\}$ to be asymptotic to one in $\Sigma$ whereas the former requires every nontrivial trajectory below $\Sigma$ to be asymptotic to one in $\Sigma$ and $\Sigma$ to attract all the points of $B\setminus\{0\}$. Obviously, the concept of a modified carrying simplex is more general and it includes carrying simplex as a particular class.

\begin{definition}
	The map $T: C\to C$ defined by (\ref{e2}) is said to be weakly retrotone in a subset $X\subset C$ if for $x, y\in X$ with $T(x) >T(y)$ and $T(x)-T(y) \in \dot{C}_I$ for some $I\subset I_N$, then $x>y$ and $x_i>y_i$ for all $i\in I$.
\end{definition}

Comparing this with the definition of retrotone given in section 1 we see that if $T$ is retrotone then it is weakly retrotone, but not vice versa.

\begin{theorem}\label{The2.1}
	Assume that $T$ defined by (\ref{e2}) with $T([0, r]) \subset [0, r]$ for some $r\in \dot{C}$ satisfies the following conditions:
	\begin{itemize}
		\item[(i)] For each $i\in I_N$, the map $T$ restricted to $X_i$ has a fixed point $q_ie_i$ with $q_i>0$ and $q\ll r$.
		\item[(ii)] $T$ is weakly retrotone and locally one to one in $[0, r]$.
		\item[(iii)]  For any $x, y\in [0, r]$, if $T(x) < T(y)$ and $T(y)-T(x) \in \dot{C}_I$ for some $I\subset I_N$ then, for each $j\in I$, either $x_j=0$ or $f_j(x) >f_j(y)$.
	\end{itemize}
	Then 0 is a repellor with the basin of repulsion $\mathcal{B}(0)\subset [0, r]$, (\ref{e1}) has a unique modified carrying simplex $\Sigma$ and $\Sigma = \overline{\mathcal{B}(0)} \setminus (\{0\}\cup\mathcal{B}(0))$. Moreover, for each $p\in\Sigma$ and every $q\in [0, r]\setminus\{0\}$ with $q<p$, we have $\alpha(q)\subset \pi_i$ provided $q_i< p_i$.
\end{theorem}

\noindent {\bf Remark 2.1} 
\begin{itemize}
	\item[(a)] Condition (i) of Theorem \ref{The2.1} is the same as that of Theorem \ref{The1.1} but conditions (ii) and (iii) are weaker than those of Theorem \ref{The1.1}. 
	\item[(b)] Condition (ii) and the definition (\ref{e2}) imply that $T: [0, r]\to T([0, r])$ is a homeomorphism. This follows from the local one to one property of $T$ on $[0, r]$, $T(x) = 0$ if and only if $x=0$, and Lemma 4.1 in \cite{Her}. 
	\item[(c)] Condition (ii) implies that, for each $i\in I_N$, the function $T_i(se_i)$ is strictly increasing for $s\in [0, r_i]$. Indeed, from (b) above we know that $T$ is one to one on $[0, r]$. As $T_j(se_i) =0$ and $T_i(se_i) > 0$ for all $j\in I_N\setminus\{i\}$ and $s\in (0, r_i]$, the one to one property of $T$ ensures that $T_i(s_1e_i) \not= T_i(s_2e_i)$ for $0< s_1<s_2\leq r_i$. By (ii) we must have  $T_i(s_1e_i) < T_i(s_2e_i)$ for $0< s_1<s_2\leq r_i$. By continuity, $T_i(se_i)$ is strictly increasing for $s\in [0, r_i]$.
	\item[(d)] Conditions (ii) and (iii) imply that, for each $i\in I_N$, $f_i(se_i)$ is strictly decreasing for $s\in [0, r_i]$. Indeed, for $0< s_1<s_2\leq r_i$, from (c) above we see that $0< T_i(s_1e_i) < T_i(s_2e_i)$ and $T(s_1e_i) < T(s_2e_i)$. From (iii) we have $f_i(s_1e_i) > f_i(s_2e_i)$. By continuity of $f$, $f_i(se_i)$ is strictly decreasing for $s\in [0, r_i]$.
	\item[(e)] The conclusion that the origin is a repellor immediately follows from conditions (i)--(iii). In fact, condition (i) implies that $f_i(q_ie_i) = 1$ for all $i\in I_N$. From (d) above we have $f_i(0) > 1$ for all $i\in I_N$. As each $f_i(0)$ is an eigenvalue of $DT(0)$, all eigenvalues of $DT(0)$ are greater than 1 so 0 is a repellor.
	\item[(f)] From the conclusion we see that for each $p\in\Sigma$, there is a nonempty $I\subset I_N$ such that $p\in \dot{C}_I$. Then, for each $q\in [0, r]$ with $q\ll_I p$, we have $\alpha(q) \subset \pi_i$ for all $i\in I$. As $q\ll_I p$ and $p\in \dot{C}_I$ imply that $q\in C_I$, we have $\alpha(q) \subset \cap^N_{i=1}\pi_i = \{0\}$, so $\alpha(q) =\{0\}\not\subset\Sigma$. Since $\alpha(q)\subset\Sigma$ if $q\in\Sigma$ by the invariance and compactness of $\Sigma$, we must have $q\not\in\Sigma$. This shows that $\Sigma$ is unordered in a strict sense: for any nonempty $I\subset I_N$, any $p\in \dot{C}_I$ and any $q\ll_I p$, we cannot have both $p\in\Sigma$ and $q\in \Sigma$. In other words, $\Sigma$ is unordered in the sense of $\ll_I$ for any nonempty $I\subset I_N$: there are no distinct points $p, q\in \Sigma\cap C_I$ such that $p\ll_I q$. However, due to the possibility of $\frac{\partial f_i(x)}{\partial x_j} = 0$ for some $i\not= j$ and some $x$, $\Sigma$ does allow ordered points on it in the sense of $<$, i.e. $p, q\in \Sigma$ with $p <q$. This is demonstrated by the trivial example below.
\end{itemize}

\noindent {\bf Example} Consider the system (\ref{e1}) with $T$ given by
\begin{equation}\label{e5}
T_i(x) = x_ig_i(x_i), i\in I_N,
\end{equation}
where each $g_i: \R_+\to\R$ is positive, continuous, $0 <g_i(u) <1$ for $u\geq r_i> q_i>0$, $g_i(q_i) = 1$, $g_i\in C^1([0, r_i], \R)$, $g'_i(u)<0$, and $g_i(u) +ug'_i(u) >0$ for $u\in [0, r_i]$. Then $T$ satisfies all the conditions of Theorem \ref{The2.1}, so it has a unique modified carrying simplex $\Sigma$. Note that system (\ref{e1}) with $T$ defined by (\ref{e5}) is a trivial case of  (\ref{e1}) with $T$ defined by (\ref{e2}) when there is no interaction between distinct component equations of the system. Since $q_i$ is the globally attracting equilibrium of the $i$th component equation on the positive $x_i$-axis, $\Sigma$ is the upper boundary surface of the cell $[0, q]$, i.e.
\[
\Sigma = \{x\in [0, q]: x_i = q_i \; \textup{for some}\; i\in I_N\}.
\]  
Clearly, $q\in\Sigma$ and for each $p\in\Sigma\setminus\{q\}$, we have $p<q$. Thus, ordered points in the sense of $<$ are permitted on $\Sigma$.

Now utilising $DT$ and $Df$, we obtain conditions which guarantee conditions (ii) and (iii) and the following version of Theorem \ref{The2.1} with easily checkable conditions. Consider the matrix $M(x)$ given by (\ref{e4}) and
\begin{equation}
\tilde{M}(x)= (\tilde{M}_{ij}(x))= \left(-\frac{x_j}{f_i(x)}\frac{\partial f_i}{\partial x_j}(x)\right)_{N\times N}. \label{E4}
\end{equation}

\begin{theorem}\label{The2.2}
	Assume that $T$ given by (\ref{e2}) satisfies the following conditions:
	\begin{itemize}
		\item[(i)] For each $i\in I_N$, the map $T$ restricted to $X_i$ has a fixed point $q_ie_i$ with $q_i>0$ and $q\ll r$ for some $r\in \dot{C}$.
		\item[(ii)] The entries of the Jacobian $Df$ satisfy
		\begin{equation}\label{e6}
		\forall x\in [0, r], \forall i, j\in I_N,\; \frac{\partial f_i}{\partial x_j}(x) \leq 0,
		\end{equation}
		and $f_i$ is strictly decreasing in $x_i\in [0, r_i]$ for $x\in [0, r]$.
		\item[(iii)] For each $x\in [0, q]\setminus\{0\}$, either $\rho(M(x))<1$ for $M(x)$ given by (\ref{e4}) or $\rho(\tilde{M}(x))<1$ for $\tilde{M}(x)$ given by (\ref{E4}).
	\end{itemize}
	Then 0 is a repellor with the basin of repulsion $\mathcal{B}(0)\subset [0, r]$, (\ref{e1}) has a unique modified carrying simplex $\Sigma$ and $\Sigma = \overline{\mathcal{B}(0)} \setminus (\{0\}\cup\mathcal{B}(0))$. Moreover, for each $p\in\Sigma$ and every $q\in [0, r]\setminus\{0\}$ with $q<p$, we have $\alpha(q)\subset \pi_i$ provided $q_i< p_i$.
\end{theorem}

\noindent {\bf Remark 2.2}
\begin{itemize}
	\item[(a)] When $\frac{\partial f_i}{\partial x_i} \leq 0$, a sufficient condition for $f_i$ to be strictly decreasing for $x_i\in [0, r_i]$, $x\in [0, r]$ with $x_j$ fixed for all $j\in I_N\setminus\{i\}$, is that the set $Z_i$ of zeros of $\frac{\partial f_i}{\partial x_i}$ in $[0, r_i]$ is either empty or finite or infinite with only a finite number of accumulation points. In particular, when each $Z_i$ is empty, condition (ii) in Theorem \ref{The2.2} can be replaced by
	\begin{itemize}
		\item[(ii)*] For all $i, j\in I_N$, the entries of the Jacobian $Df$ satisfy
		\begin{equation}\label{e6''}
		\forall x\in [0, r], \; \frac{\partial f_i}{\partial x_i}(x) < 0, \;\frac{\partial f_i}{\partial x_j}(x) \leq 0.
		\end{equation}
	\end{itemize} 
	\item[(b)] Comparing Theorem \ref{The2.2} with Theorem \ref{The1.2}, we see that condition (i) of Theorem \ref{The2.2} is the same as (i) of 
	Theorem \ref{The1.2} and (iii) of Theorem \ref{The2.2} has one more choice than (iii) of Theorem \ref{The1.2}, but condition (ii) of Theorem \ref{The2.2} only requires each entry of $Df$ to be nonnegative instead of $N^2$ entries of $Df$ to be strictly negative in Theorem \ref{The1.2}, plus the strictly decreasing requirement of each $f_i$ in $x_i$. Even if (ii) is replaced by the stronger condition (ii)* above, it only requires $N$ diagonal entries of $Df$ to be negative.  From this point of view, with a trade off of having a modified carrying simplex rather than the well known carrying simplex, we have significantly reduced the cost and generalised the existing results.
	\item[(c)] Under condition (ii) of Theorem \ref{The2.2}, if
	\begin{equation}\label{E7'}
	f_i(x) + \sum_{j=1}^Nx_i\frac{\partial f_i}{\partial x_j}(x)>0 \quad \forall i\in I_N,
	\end{equation}
	using one type of matrix norm we have
	\[
	\|M(x)\|=\max_{i\in I_N}\sum^N_{j=1}\left|\frac{x_i}{f_i(x)}\frac{\partial f_i(x)}{\partial x_j}\right| <1.
	\] 
	By Theorem 6.1.3 in \cite{Lan}, $\rho(M(x))\leq \|M(x)\|$. Thus, (\ref{E7'}) is a sufficient condition for $\rho(M(x))<1$. By the same reason, if
	\begin{equation}\label{E7''}
	f_i(x) + \sum_{j=1}^Nx_j\frac{\partial f_i}{\partial x_j}(x)>0 \quad \forall i\in I_N,
	\end{equation}
	then $\rho(\tilde{M}(x))\leq \|\tilde{M}(x)\| <1$. Therefore, condition (iii) of Theorem \ref{The2.2} is met if (\ref{E7'}) or (\ref{E7''}) holds for each $x\in [0, q]\setminus\{0\}$.
\end{itemize}

\begin{corollary}\label{Cor2.1}
	Under the conditions of Theorem \ref{The2.1} or Theorem \ref{The2.2}, the following conclusions hold.
	\begin{itemize}
		\item[(i)] For any periodic orbit $\gamma\subset\Sigma$, the points on $\gamma$ are unordered, i.e. if $p, q\in \gamma$ with $p\leq q$ then $p=q$.
		\item[(ii)] For any $x\in\Sigma$, if there are two points $p, q\in \gamma(x)$ satisfying $p<q$ then $\alpha(x)$ consists of either a single fixed point or a periodic orbit.
	\end{itemize}
\end{corollary}

\begin{proof}
	(i) Suppose there are two points $p, q\in \gamma$ satisfying $p<q$. Then there is at least one $i\in I_N$ such that $p_i<q_i$. From Theorem \ref{The2.1} we have $\alpha(p)\subset \pi_i$, so $q\not\in \alpha(p)$, a contradiction to $q\in\gamma = \alpha(p)$ due to the periodicity of $\gamma$. Therefore, $\gamma$ is unordered.
	
	(ii) By $x\in\Sigma$ we have $\gamma(x)\subset \Sigma$ and $\alpha(x)\subset\Sigma$. Since $p, q\in \gamma(x)$ with $p<q$, we have $T(T^{-1}(p)) = p < q = T(T^{-1}(q))$. Then the weakly retrotone property of $T$ implies that $T^{-1}(p)<T^{-1}(q)$ and $T^{-n}(p) < T^{-n}(q)$ for all $n\in\N$. For each $i\in I_N$, if there is an $n\in\N$ such that $(T^{-n}(p))_i < (T^{-n}(q))_i$, by Theorem \ref{The2.1} we have $\alpha(x) = \alpha(p) \subset \pi_i$; otherwise, we have $(T^{-n}(p))_i = (T^{-n}(q))_i$ for all $n\in\N$. Thus, there is a proper subset $I\subset I_N$ such that $\alpha(x)\subset \pi_i$ for each $i\in I$ and $(T^{-n}(p))_j = (T^{-n}(q))_j$ for all $n\in\N$ and $j\in I_N\setminus I$. As $p$ and $q$ are two distinct
	points on $\gamma(x)$, there is an $n_1>0$ such that either $T^{n_1}(p) = q$ or $T^{n_1}(q) = p$. Hence, since the component $(T^n(p))_j$ is an $n_1$-periodic function for $n\in\N$ for each $j\in I_N\setminus I$, we obtain
	\begin{eqnarray*}
		\alpha(x) &=& \{T^k(y): k\in \{0, 1, \ldots, n_1-1\}, (T^k(y))_i=0, i\in I;\\
		& & (T^k(y))_j=(T^{k-n_1}(p))_j, j\in I_N\setminus I\}
	\end{eqnarray*}
	Therefore, $\alpha(x)$ consists of either a single fixed point or a periodic orbit. 
\end{proof}

\noindent {\bf Remark 2.3} Just as we mentioned after Theorem \ref{The1.1}, Theorems \ref{The2.1} and \ref{The2.2} can be only applied to systems on the space $[0, r]\subset C$ if no condition for $T$ on $C\setminus[0, r]$ is given. However, a simple additional condition
\begin{equation}\label{e6'}
\forall i\in I_N, \forall x\in C \; \textup{with}\; x_i\geq r_i, 0 < f_i(x) <1
\end{equation} 
guarantees that for each compact set $S\subset C$ there is a $k\in\N$ such that $T^k(S)\subset [0, r]$, so that Theorems \ref{The2.1} and \ref{The2.2} can be applied directly to systems on $C$.

In general, for any topological space $X$, a system $x(n) = F^n(x)$ for $x\in X, n\in \N$ with a map $F: X\to X$, and a compact invariant set $A\subset X$, $A$ is called a global attractor of the system if $A$ attracts the points of any bounded set $B\subset X$ uniformly. For our system (\ref{e1}) with (\ref{e2}) on $C$ under the conditions of Theorem \ref{The2.1} or Theorem \ref{The2.2}, since $0$ is a repelling fixed point, by saying that \emph{$\Sigma$ is a global attractor of the system in} $[0, r]\setminus\{0\}$ ($C\setminus\{0\}$), we mean $\Sigma$ uniformly attracts the points of any bounded set $B\subset [0, r]\setminus\{0\}$ ($B\subset C\setminus\{0\}$) that is bounded away from $0$, i.e. $\overline{B}\subset [0, r]\setminus\{0\}$ ($\overline{B}\subset C\setminus\{0\}$).

\begin{corollary}\label{Cor2.2}
	Under the conditions of Theorem \ref{The2.1} or Theorem \ref{The2.2}, the modified carrying simplex $\Sigma$ is a global attractor in $[0, r]\setminus\{0\}$. In addition, if (\ref{e6'}) holds, then $\Sigma$ is a global attractor of the system in $C\setminus\{0\}$.
\end{corollary}

Before we prove Theorem \ref{The2.1}, Theorem \ref{The2.2} and Corollary \ref{Cor2.2} in section 5, we present an application of Theorem \ref{The2.2} in next section.

\section{Criteria for dominance and vanishing species}\label{Sec3}
In this section, we consider (\ref{e2}) and give sufficient conditions for dominance of some species under the assumption that the conditions of Theorem \ref{The2.2} are met.

Viewing (\ref{e2}) as a population model for $N$ competitive species, we say that the $j$th species is {\it dominated} or {\it vanishing} if $\lim_{n\to +\infty}x_j(n) = 0$ for all $x\in \dot{C}$; we say that the $j$th species is {\it dominant} if $\liminf_{n\to\infty}x_j(n)>\delta>0$ for all $x\in\dot{C}$ and all other species are vanishing. Let
\begin{equation}\label{e7}
\Gamma_i = \{x\in C: f_i(x) =1\}, i\in I_N.
\end{equation}
Under the general assumptions for (\ref{e2}), each $\Gamma_i$ is a closed set and an $(N-1)$-dimensional hypersurface. In this section, we assume that each $\Gamma_i$ is simply connected and divides $C$ into three mutually exclusive subsets $\Gamma^+_i, \Gamma_i$ and $\Gamma^-_i$ with $0\in \Gamma^-_i$. Then the closure of $\Gamma^-_i$ is $\overline{\Gamma^-_i} = \Gamma^-_i\cup \Gamma_i$. But if we consider the restriction of $\Gamma_i$ to $[0, r]$, this assumption is met if the conditions of Theorem \ref{The2.2} hold: each $\Gamma_i\cap [0, r]$ is a simply connected closed set and an $(N-1)$-dimensional hypersurface such that $\Gamma_i^-\cap [0, r]$ is strictly below $\Gamma_i$ and $\Gamma_i^+\cap [0, r]$ is strictly above $\Gamma_i$. 

Under the conditions of Theorem \ref{The2.2}, let $Q_i = q_ie_i$, the fixed point of $T$ on $X_i$, the positive half $x_i$-axis.

\begin{theorem}\label{The3.1}
	Assume that (\ref{e6'}) and the conditions of Theorem \ref{The2.2} hold. 
	\begin{itemize}
		\item[(a)] If for some $i\in I_N$ and all $j\in I_N\setminus\{i\}$, $\frac{\partial f_i}{\partial x_i}(Q_i) <0$ and $\Gamma_i\cap [0, r]$ is strictly below $\Gamma_j$, then $\lim_{n\to +\infty}x_i(n) = 0$ for all $x\in C\setminus X_i$ so the $i$th species is vanishing.
		\item[(b)] If for some $i\in I_N$ and all $j\in I_N\setminus\{i\}$, $\frac{\partial f_i}{\partial x_i}(Q_i) <0$ and $\Gamma_i\cap [0, r]$ is strictly above $\Gamma_j$, then the $i$th species is dominant and the axial fixed point $Q_i$ is globally asymptotically stable.
	\end{itemize}
\end{theorem}

\begin{proof}
	By Theorem \ref{The2.2} the system has a unique modified carrying simplex $\Sigma$.
	
	(a) Under the assumption that $\Gamma_i\cap [0, r]$ is strictly below $\Gamma_j$ for all $j\in I_N\setminus\{i\}$, we first claim that 
	\begin{equation}\label{e8}
	\Gamma^-_i \cap \Sigma = \emptyset
	\end{equation}
	so that $\Gamma^-_i\cap [0, r]$ is strictly below $\Sigma$ and $\Sigma$ is above $\Gamma_i$. Indeed, if (\ref{e8}) were not true then we would have a point $p\in (\Gamma^-_i\cap\Sigma$). As $0\not\in \Sigma$, we have $p\not= 0$ and a nonempty $J\subset I_N$ such that $p_j>0$ if and only if $j\in J$. Since $p$ is below $\Gamma_i$ and $\Gamma_i\cap [0, r]$ is strictly below $\Gamma_j$ for all $j\in I_N\setminus\{i\}$, $p$ is below $\Gamma_j$ for all $j\in I_N$. Let $u = T(p)$. Then 
	\[
	\forall j\in J, u_j = T_j(p) = p_jf_j(p) > p_j; \forall k\in I_N\setminus J, u_k = p_k = 0,
	\]
	so $p\ll_J u$. As $\Sigma$ is invariant and $p\in \Sigma$, we have $u\in \Sigma$. Then, by Theorem \ref{The2.2}, we have $\alpha(p) = \{0\}$ so $p\in \mathcal{B}(0)$, a contradiction to $p\in \Sigma = \overline{\mathcal{B}(0)}\setminus (\{0\}\cup\mathcal{B}(0))$. This shows our claim (\ref{e8}).
	
	Since the axial fixed point $Q_i$ is below $\Gamma_j$ for all $j\in I_N\setminus\{i\}$, the Jacobian matrix $DT(Q_i)$ has $N-1$ eigenvalues $f_j(Q_i) >1$ for $j\in I_N\setminus\{i\}$ and one eigenvalue $1+q_i\frac{\partial f_i}{\partial x_i}(Q_i)$. By assumption, (\ref{e4}) and (\ref{E4}), the only nonzero eigenvalue of $M(Q_i)$ and $\tilde{M}(Q_i)$ is $-q_i\frac{\partial f_i}{\partial x_i}(Q_i)>0$, so $\rho(M(Q_i))=\rho(\tilde{M}(Q_i))=-q_i\frac{\partial f_i}{\partial x_i}(Q_i)$. By condition (iii) of Theorem \ref{The2.2}, we have $0 <1+q_i\frac{\partial f_i}{\partial x_i}(Q_i) < 1$.  So $Q_i$ is a saddle point in $C$ with $X_i$ as its one-dimensional stable manifold and a repellor on $\Sigma$. Thus, to show that $\lim_{n\to +\infty}x_i(n) = 0$ for all $x\in C\setminus X_i$, by the definition of modified carrying simplex, we need only show that $\lim_{n\to +\infty}x_i(n) = 0$ for all $x\in C\setminus X_i$ on or above $\Sigma$, i.e. $x\in (\Sigma\cup\Sigma^+)\setminus X_i$. 
	
	Now for any $x\in C$ with $x_i >r_i$, the assumption (\ref{e6'}) ensures that $x(n)\in [0, r]$ for large enough $n\in\N$. Without loss of generality, we only consider $x\in (\Sigma\cup(\Sigma^+\cap [0, r]))\setminus X_i$. We first show that the set $(\Sigma\cup(\Sigma^+\cap [0, r]))\setminus X_i$ is positively invariant. From the proof of Theorem \ref{The2.2} given in section 5 we shall see that the conditions of Theorem \ref{The2.2} imply the conditions of Theorem \ref{The2.1}. Thus, $[0, r]$ is positively invariant and, by Remark 2.3 (b), $T: [0, r] \to T([0, r])$ is a homeomorphism. As $0$ is a repellor with the basin of repulsion $\mathcal{B}(0)\subset [0, r]$, we shall see in section 5 (Lemma \ref{Lem5.4}) that $\overline{\mathcal{B}(0)}$ is invariant. Thus, $T$ maps the set
	\[
	[0, r]\setminus \overline{\mathcal{B}(0)} = [0, r] \setminus (\Sigma\cup\mathcal{B}(0)\cup\{0\}) = [0, r]\setminus(\Sigma\cup\Sigma^-) = [0, r]\cap\Sigma^+
	\]
	into itself. As $\Sigma$ is invariant, $\Sigma\cup(\Sigma^+\cap[0, r])$ is positively invariant. For each $x\in \Sigma\cup(\Sigma^+\cap[0, r])\setminus X_i$, there is a $j\in I_N \setminus\{i\}$ such that $x_j>0$, so $T_j(x) = x_jf_j(x) >0$. Thus, $T(x)\in \Sigma\cup(\Sigma^+\cap[0, r])\setminus X_i$. This shows the positive invariance of $\Sigma\cup(\Sigma^+\cap[0, r])\setminus X_i$. 
	
	By (\ref{e8}), $\Sigma$ is above $\Gamma_i$, so $\Sigma\cup(\Sigma^+\cap[0, r])\setminus X_i$ is above $\Gamma_i$. Thus, for $x\in (\Sigma\cup(\Sigma^+\cap [0, r]))\setminus X_i$, $x(n) = T^n(x)\in (\Sigma\cup(\Sigma^+\cap [0, r]))\setminus X_i$, so $x(n)$ is on or above $\Gamma_i$ for all $n\in\N$. Hence,
	\[
	\forall n\in\N, x_i(n+1) = T_i(x(n)) = x_i(n)f_i(x(n)) \leq x_i(n).
	\]
	This shows that $\{x_i(n)\}$ is a bounded monotone nonincreasing sequence, so there is an $x_0\geq 0$ such that $\lim_{n\to +\infty}x_i(n) = x_0$. Suppose $x_0>0$. Then, for each $y\in \omega(x)\subset\Sigma$, we have $y_i =x_0$. As $T^n(y) \in \omega(x)$ for all integer $n$, we have $T_i(y) = y_if_i(y) =x_0 = y_i$ so $f_i(y) =1$ and $y\in \Gamma_i$. Therefore, $\omega(x)\subset \Gamma_i\cap\Sigma$. If $\omega(x) =\{Q_i\}$, as $Q_i$ is below $\Gamma_j$ for all $j\in I_N\setminus\{i\}$, there is a $\delta>0$ such that the closure $\overline{O(Q_i, \delta)}\cap [0, r]$ of the open ball centred at $Q_i$ with radius $\delta$ restricted to $[0, r]$, i.e. $O(Q_i, \delta)\cap [0, r]$, is strictly below $\Gamma_j$ for all $j\in I_N\setminus\{i\}$. Let
	\[
	m_0 = \min\{f_j(u): u\in \overline{O(Q_i, \delta)}\cap [0, r], j\in I_N\setminus\{i\}\}.
	\]
	Then $m_0 >1$. Since $\lim_{n\to +\infty}x(n) = Q_i$, there is $n_1\in\N$ such that $x(n) \in O(Q_i, \delta)\cap [0, r]$ for $n\geq n_1$. As $x\not\in X_i$, we have $x_j>0$ for some $j\in I_N\setminus\{i\}$. Then, for this $j$ and all $n\geq 1$,
	\[
	x_j(n+n_1) = T_j(x(n-1+n_1) = x_j(n-1+n_1)f_j(x(n-1+n_1) \geq m_0 x_j(n-1+n_1),
	\]
	so $x_j(n+n_1) \geq m^n_0 x_j(n_1) \to +\infty$ as $n\to +\infty$, a contradiction to the boundedness of $\{x(n)\}$. This contradiction shows the existence of a point $y\in\omega(x) \setminus \{Q_i\}$. Since $Q_i$ is the unique intersection point of $\Sigma$ with $X_i$ and $y\in \Sigma\setminus\{Q_i\}$, we have $y\not\in X_i$ so $y_j>0$ for some $j\in I_N\setminus\{i\}$. Since $\Gamma_i\cap [0, r]$ is strictly below $\Gamma_j$ and $\omega(x)\subset \Gamma_i\cap\Sigma \subset \Gamma_i\cap [0, r]$, $\omega(x)$ is strictly below $\Gamma_j$. Let
	\[
	\rho = \min_{u\in\omega(x)}f_j(u).
	\]
	Then, by the continuity of $f$ and the compactness of $\omega(x)$, $\rho >1$ and 
	\[
	y_j(n+1) = T_j(y(n)) = y_j(n)f_j(y(n))\geq \rho y_j(n). 
	\]
	Thus,
	\[
	y_j(n)\geq \rho^n y_j \to +\infty \;(n\to +\infty),
	\]
	a contradiction to the boundedness of $\omega(x)$. This contradiction shows that we must have $x_0 = 0$, i.e. $\lim_{n\to +\infty}x_i(n) = 0$ for all $x\in C\setminus X_i$.
	
	(b) Under the condition that $\Gamma_i\cap [0, r]$ is strictly above $\Gamma_j$ for every $j\in I_N\setminus\{i\}$, we first show that $\Sigma$ is below $\Gamma_i$ by assuming the opposite: there is a point $p\in \Sigma\cap\Gamma^+_i$. As $p\not=0$, there is a nonempty $J\subset I_N$ as the support of $p$. As $p$ is above $\Gamma_i$ and $\Gamma_i\cap [0, r]$ is strictly above $\Gamma_j$ for all $j\in I_N\setminus\{i\}$, $p$ is above $\Gamma_j$ for all $j\in I_N$. Thus,
	\[
	\forall j\in J, T_j(p) = p_jf_j(p) < p_j,
	\]
	so $T(p)\ll_J p$. By Theorem \ref{The2.2}, $\alpha(T(p)) =\{0\}$ so $T(p) \in \mathcal{B}(0)$, a contradiction to $T(p)\in \Sigma = \overline{\mathcal{B}(0)}\setminus(\{0\}\cup\mathcal{B}(0))$. This shows that $\Sigma$ must be below $\Gamma_i$.
	
	We need only show that $Q_i$ is stable and attracts all the points of $C\setminus\pi_i$ as the dominance of the $i$th species is implied by the global attraction of $Q_i$. As $Q_i$ is above $\Gamma_j$ for all $j\in I_N\setminus\{i\}$, we have $f_j(Q_i) \in (0, 1)$ for all $j\in I_N\setminus\{i\}$. By the assumption $\frac{\partial f_i}{\partial x_i}(Q_i)<0$ and condition (iii) of Theorem \ref{The2.2}, $1+q_i\frac{\partial f_i}{\partial x_i}(Q_i) \in (0, 1)$. Thus, every eigenvalue of $DT(Q_i)$ is in $(0, 1)$ so $Q_i$ is asymptotically stable. To show the global attraction of $Q_i$ in $C\setminus\pi_i$, by the assumption (\ref{e6'}) and the definition of a modified carrying simplex, we need only show that $\lim_{n\to +\infty}x(n) =Q_i$ for all $x\in (\Sigma\cup(\Sigma^+\cap [0, r]))\setminus\pi_i$. 
	
	If $x\in\Sigma\setminus\pi_i$, as $x(n)\in \Sigma$ for all $n\in\N$ and $\Sigma$ is below $\Gamma_i$, the sequence $\{x_i(n)\}$ is bounded and monotone nondecreasing. Thus, there is a $\beta>0$ such that $\lim_{n\to +\infty}x_i(n) = \beta$. For any $y\in \omega(x)$, we have $T_i(y) = y_if_i(y) =\beta = y_i$ so $f_i(y) = 1$. Thus, $y\in\Gamma_i$ and $\omega(x)\subset \Gamma_i$. We claim that $\omega(x) = \{Q_i\}$. To verify this claim, as $\omega(x)$ is compact, $\omega(x)\subset \Gamma_i$ and $\Gamma_i\cap [0, r]$ is strictly above $\Gamma_j$ for all $j\in I_N\setminus\{i\}$, there is a $\delta >0$ such that the closure $\overline{O(\omega(x), \delta)}\cap [0, r]$ of the open set $O(\omega(x), \delta)\cap [0, r]$ with
	\[
	O(\omega(x), \delta) = \{z\in C: \|z-u\| < \delta \; \textup{for some}\; u\in\omega(x)\}
	\]
	is strictly above $\Gamma_j$ for all $j\in I_N\setminus\{i\}$. Let
	\[
	\mu = \max\{f_j(u): u\in\overline{O(\omega(x), \delta)}\cap [0, r], j\in I_N\setminus\{i\}\}.
	\]
	By the continuity of $f$ and the compactness of $\overline{O(\omega(x), \delta)}\cap [0, r]$, we have $0<\mu <1$. By the definition of $\omega(x)$, there is an integer $N_1\geq 0$ such that $x(n) \in O(\omega(x), \delta)\cap [0, r]$ for all $n> N_1$. Let $J\subset I_N$ be the support of $x$. Then $x_j(n) >0$ for all $n\in\N$ if and only if $j\in J$. If $J = \{i\}$ then $x=Q_i$ and the above claim is obviously true. Now suppose $\{i\}\subset J\not= \{i\}$. Then,
	\[
	\forall n>N_1, \forall j\in J\setminus\{i\}, x_j(n+1) = x_j(n)f_j(x(n)) \leq \mu x_j(n),
	\]
	so $x_j(n+N_1) \leq \mu^n x_j(N_1)\to 0$ as $n\to +\infty$. This shows that $\omega(x) = \Sigma\cap X_i=\{Q_i\}$.
	
	Now suppose $x\in (\Sigma^+\cap [0, r])\setminus\pi_i$. By the asymptotic stability of $Q_i$ there is a $\delta>0$ such that every point $z\in O(Q_i, \delta)\cap C$ is attracted to $Q_i$, i.e. $\lim_{n\to +\infty}z(n) = Q_i$. Thus, as long as $Q_i\in\omega(x)$, there is an $m\in\N$ such that $x(m)\in O(Q_i, \delta)\cap C$ so that $\omega(x) = \{Q_i\}$. We now prove $Q_i\in\omega(x)$ by contradiction. If $Q_i\not\in\omega(x)$ and there is a $y\in\omega(x)\subset \Sigma$ with $y_i>0$, by the previous paragraph we have $\lim_{n\to +\infty} y(n) = Q_i$. As $\omega(x)$ is compact and $y(n)\in\omega(x)$ for all $n\in\N$, we have $Q_i\in \omega(x)$, which contradicts the condition $Q_i\not\in\omega(x)$. Thus, if $Q_i\not\in\omega(x)$ then $\omega(x)\subset \Sigma\cap\pi_i$. If $\omega(x)$ is strictly below $\Gamma_i$, by definition of $\omega(x)$ there is a $K\in\N$ such that $x(n)$ is below $\Gamma_i$ for all $n\geq K$. Thus, $\{x_i(n)\}$ is an increasing sequence for $n\geq K$ so that $x_i(n)\geq x_i(K)>0$ for $n\geq K$ and each $y\in\omega(x)$ satisfies $y_i\geq x_i(K)>0$, a contradiction to the assumption $\omega(x)\subset \pi_i$. If there is a point $p\in \omega(x)\subset (\Sigma\cap\pi_i)$ on or above $\Gamma_i$, there is a nonempty $J\subset I_N\setminus\{i\}$ as the support of $p$. As $\Gamma_i\cap [0, r]$ is strictly above $\Gamma_j$ for all $j\in I_N\setminus\{i\}$, $p$ is above $\Gamma_j$ for all $j\in I_N\setminus\{i\}$ so $T(p) \ll_J p$. This leads us to $\alpha(T(p)) = \{0\}$ by Theorem \ref{The2.2}, so $T(p) \in \mathcal{B}(0)$, a contradiction to $T(p)\in\Sigma = \overline{\mathcal{B}(0)}\setminus(\{0\}\cup\mathcal{B}(0))$. These contradictions show that we must have $Q_i\in\omega(x)$ so $\omega(x) = \{Q_i\}$.  
\end{proof}

We note that Theorem 2.3 in \cite{Her} is consistent with our Theorem \ref{The3.1} (b) but under the stronger conditions of Theorem \ref{The1.1}. While Theorem \ref{The3.1} used one surface $\Gamma_i\cap [0, r]$ comparing with the other $N-1$ surfaces $\Gamma_j$ to obtain one species vanishing, our next result repeat such a condition several times to get multiple species vanishing.

\begin{theorem}\label{The3.2}
	Assume that (\ref{e6'}) and the conditions of Theorem \ref{The2.2} hold. Assume also the existence of an integer $k\in I_N\setminus\{N\}$ such that  for all $i\in\{1, \ldots, k\}$, $\frac{\partial f_i}{\partial x_i}(Q_i)<0$ and
	\begin{equation}\label{e9}
	\forall j\in\{i+1, \ldots, N\}, (\cap^{i-1}_{\ell=1}\pi_{\ell})\cap \Gamma_i\cap[0, r]\; \textit{is strictly below}\; \Gamma_j.
	\end{equation}
	Then the $i$th species is dominated for all $i\in \{1, \ldots, k\}$. In addition, if $k = N-1$ and $\frac{\partial f_N}{\partial x_N}(Q_N) <0$, then the $N$th species is dominant and the $N$th axial fixed point $Q_N$ is globally asymptotically stable in $C$.
\end{theorem}

\noindent {\bf Remark 3.1} Here the symbol $\cap_{i\in\emptyset}\pi_i$ is deemed as $C$. So, (\ref{e9}) for $i=1$ is simplified as
\begin{equation}\label{e10}
\forall j\in\{2, \ldots, N\}, \Gamma_1\cap [0, r]\; \textup{is strictly below}\; \Gamma_j.
\end{equation}	

\begin{proof}[Proof of Theorem \ref{The3.2}]
	For $k>1$ we first prove that
	\begin{equation}\label{e13}
	\forall i\in\{2, \ldots, k\}, (\cap_{\ell=1}^{i-1}\pi_{\ell})\cap\Gamma^-_i \cap\Sigma = \emptyset.
	\end{equation}
	The proof of (\ref{e13}) is similar to that of (\ref{e8}). Suppose (\ref{e13}) is not true. Then, for some $i\in \{2, \ldots, k\}$, there exists a point $u\in (\cap_{\ell=1}^{i-1}\pi_{\ell})\cap\Gamma^-_i \cap\Sigma$, so $u$ is below $\Gamma_i$. By (\ref{e9}), $(\cap_{\ell=1}^{i-1}\pi_{\ell})\cap [0, r]\cap\Gamma_i$ is strictly below $\Gamma_j$ for all $j\in\{i+1, \ldots, N\}$. Thus, $u$ is below $\Gamma_j$ for all $j\in\{i, \ldots, N\}$. Note that $u\in(\cap_{\ell=1}^{i-1}\pi_{\ell})\cap\Sigma$ implies $u>0$ and $u_j=0$ for all $j\in\{1, \ldots, i-1\}$. Thus, there is a nonempty $J\subset \{i, \ldots, N\}$ such that $u_j>0$ if and only if $j\in J$. Then, 
	\[
	\forall j\in J, T_j(u) = u_jf_j(u) >u_j,
	\]
	so $u\ll_J T(u)$. As $T(u)\in\Sigma$, by Theorem \ref{The2.2} we obtain $u\in\mathcal{B}(0)$, a contradiction to $u\in \Sigma = \overline{\mathcal{B}(0)}\setminus(\{0\}\cup\mathcal{B}(0))$. This contradiction shows the truth of (\ref{e13}).
	
	For $k\geq 1$ with $i=1$, from Remark 3.1 and Theorem \ref{The3.1} we know that the first species is vanishing, i.e. $\omega(x) \subset \Sigma\cap\pi_1$ for all $x\in C\setminus \overline{X_1}$. 
	
	Then for $k\geq 2$ with $i=2$, from (\ref{e9}) we see that $\Gamma_2\cap\pi_1\cap [0, r]$ is strictly below $\Gamma_j$ for all $j\in \{3, \ldots, N\}$. As $\omega(x) \subset \Sigma\cap\pi_1$ for all $x\in C\setminus \overline{X_1}$, we can prove that $\omega(x) \subset \Sigma\cap(\pi_1\cap\pi_2)$ for all $x\in C\setminus (\overline{X_1}\times \overline{X_2})$ (the proof is included in the general case below).
	
	Then for $k\geq 3$ with $i=3$, from (\ref{e9}) we see that $\Gamma_3\cap(\pi_1\cap\pi_2)\cap [0, r]$ is strictly below $\Gamma_j$ for all $j\in \{4, \ldots, N\}$. As $\omega(x) \subset \Sigma\cap(\pi_1\cap\pi_2)$ for all $x\in C\setminus (\overline{X_1}\times \overline{X_2})$, we can prove that $\omega(x) \subset \Sigma\cap(\cap_{\ell=1}^3\pi_{\ell})$ for all $x\in C\setminus (\overline{X_1}\times \overline{X_2}\times \overline{X_3})$ (the proof is included in the general case below).
	
	In general, for $k >1$ with $i=k$, from (\ref{e9}) we see that $\Gamma_k\cap(\cap_{\ell=1}^{k-1}\pi_{\ell})\cap [0, r]$ is strictly below $\Gamma_j$ for all $j\in\{k+1, \ldots, N\}$. Suppose we know that $\omega(x) \subset \Sigma\cap(\cap_{\ell=1}^{k-1}\pi_{\ell})$ for all $x\in C\setminus(\overline{X_1}\times\cdots\times \overline{X_{k-1}})$. We need to prove that
	\begin{equation}\label{e14}
	\forall x\in C\setminus(\overline{X_1}\times\cdots\times \overline{X_k}), \omega(x)\subset \Sigma\cap(\cap_{\ell=1}^k\pi_{\ell}).
	\end{equation}
	From condition (\ref{e6'}), Theorem \ref{The2.2} and the definition of a modified carrying simplex, instead of (\ref{e14}) we need only prove that
	\begin{equation}\label{e15}
	\forall x\in ([0, r]\setminus\Sigma^-)\setminus(\overline{X_1}\times\cdots\times \overline{X_k}), \omega(x)\subset \Sigma\cap(\cap_{\ell=1}^k\pi_{\ell}).
	\end{equation}
	The proof of (\ref{e15}) is divided into the following two steps.
	
	Step 1. We show that
	\begin{equation}\label{e16}
	\forall x\in (\cap_{\ell=1}^{k-1}\pi_{\ell})\cap([0, r]\setminus\Sigma^-)\setminus X_k, \omega(x)\subset \Sigma\cap(\cap_{\ell=1}^k\pi_{\ell}).
	\end{equation}
	From (\ref{e13}) we know that $(\cap_{\ell=1}^{k-1}\pi_{\ell})\cap\Sigma$ is above $\Gamma_k$, so $(\cap_{\ell=1}^{k-1}\pi_{\ell})\cap([0, r]\setminus\Sigma^-)\setminus X_k$ is above $\Gamma_k$. Note that $\pi_{\ell}$ and $C\setminus \pi_{\ell}$ are positively invariant for any $\ell\in I_N$. Thus, for any $x\in (\cap_{\ell=1}^{k-1}\pi_{\ell})\cap([0, r]\setminus\Sigma^-)\setminus X_k$, if $x_k= 0$ then $x\in (\cap_{\ell=1}^k\pi_{\ell})$, so $x(n)\in(\cap_{\ell=1}^k\pi_{\ell})$ for all $n\in\N$ and $\omega(x)\subset \Sigma\cap(\cap_{\ell=1}^k\pi_{\ell})$. If $x_k>0$, as both  $\cap_{\ell=1}^{k-1}\pi_{\ell}$ and  $([0, r]\setminus\Sigma^-)\setminus X_k = \Sigma\cup(\Sigma^+\cap[0, r])\setminus X_k$ (from the proof of Theorem \ref{The3.1}) are positively invariant and  $(\cap_{\ell=1}^{k-1}\pi_{\ell})\cap([0, r]\setminus\Sigma^-)\setminus X_k = (\cap_{\ell=1}^{k-1}\pi_{\ell})\cap (([0, r]\setminus\Sigma^-)\setminus X_k) $, we have $x(n)\in (\cap_{\ell=1}^{k-1}\pi_{\ell})\cap([0, r]\setminus\Sigma^-)\setminus X_k$ for all $n\in\N$, so each $x(n)$ is on or above $\Gamma_k$ for all $n\in\N$. Hence, 
	\[
	\forall n\in\N, x_k(n+1) = T_k(x(n)) = x_k(n)f_k(x(n)) \leq x_k(n).
	\]
	As $\{x_k(n)\}$ is a positive monotone nonincreasing sequence, there is a $\mu\geq 0$ such that $\lim_{n\to +\infty}x_k(n) = \mu$. Suppose $\mu>0$. Then, for each $y\in\omega(x)\subset \Sigma$ we have $y_k = \mu$. As $T(y) \in\omega(x)$, we have $\mu = T_k(y)= y_kf_k(y)=\mu f_k(y)$ so $f_k(y) = 1$ and $y\in\Gamma_k$. Therefore, $\omega(x)\subset\Gamma_k$. By the positive invariance of $(\cap_{\ell=1}^{k-1}\pi_{\ell})\cap [0, r]$, $\omega(x)\subset \Gamma_k\cap(\cap_{\ell=1}^{k-1}\pi_{\ell})\cap [0, r]$. Since the set $\Gamma_k\cap(\cap_{\ell=1}^{k-1}\pi_{\ell})\cap [0, r]$ is strictly below $\Gamma_j$ for all $j\in\{k+1, \ldots, N\}$, $\omega(x)$ is strictly below $\Gamma_j$ for all $j\in\{k+1, \ldots, N\}$. As $\omega(x)$ is compact, there is a $\delta >0$ such that $\overline{O(\omega(x), \delta)}\cap [0, r]$ is strictly below  $\Gamma_j$ for all $j\in\{k+1, \ldots, N\}$. For this $\delta$, there is an $m\in\N$ such that $x(n)\in O(\omega(x), \delta)\cap [0, r]$ for all $n\geq m$. Note that $x\in (\cap_{\ell=1}^{k-1}\pi_{\ell})\cap([0, r]\setminus\Sigma^-)\setminus X_k$ implies $x_j>0$ for some $j\in\{k+1, \ldots, N\}$. For this $j$, let
	\[
	\eta = \min\{f_j(u): u\in \overline{O(\omega(x), \delta)}\cap [0, r]\}.
	\]
	Then $\eta > 1$ and for all $n\geq 1$, 
	\[
	x_j(n+m) = T_j(x(n+m-1)) = x_j(n+m-1)f_j(x(n+m-1)) \geq \eta x_j(n+m-1).
	\]
	It follows from this that $x_j(n+m) \geq \eta^n x_j(m)\to +\infty$ as $n\to +\infty$, a contradiction to the boundedness of $\{x(n)\}$. This contradiction shows that $\mu =0$ and (\ref{e16}) follows.
	
	Step 2. Now we prove (\ref{e15}). For $x\in ([0, r]\setminus\Sigma^-)\setminus(\overline{X_1}\times\cdots\times \overline{X_k})$, we show that $\omega(x) \subset \Sigma\cap(\cap_{\ell=1}^k\pi_{\ell})$. From the supposition we know that $\omega(x) \subset \Sigma\cap(\cap_{\ell=1}^{k-1}\pi_{\ell})$. Suppose $\omega(x) \not\subset \Sigma\cap(\cap_{\ell=1}^k\pi_{\ell})$. Then either $\omega(x)=\{Q_k\}$ or there is a $y\in \omega(x)\setminus \{Q_k\}$ with $y_k > 0$.  
	
	In the former case, $\lim_{n\to +\infty}x(n) = Q_k$. Since $Q_k$ is below $\Gamma_j$ for all $j\in\{k+1, \ldots, N\}$, there is an $\varepsilon>0$ such that $\overline{O(Q_k, \varepsilon)}\cap [0, r]$ is strictly below $\Gamma_j$ for all $j\in\{k+1, \ldots, N\}$. For this $\varepsilon >0$, there is an $m\in\N$ such that $x(n)\in O(Q_k, \varepsilon)\cap [0, r]$ for all $n\geq m$. That $x\not\in \overline{X_1}\times\cdots\times \overline{X_k}$ ensures the existence of some $j\in\{k+1, \ldots, N\}$ with $x_j>0$ so that $x_j(n)>0$ for all $n\in\N$. For this $j$, let
	\[
	\eta_0 = \min\{f_j(u): u\in \overline{O(Q_k, \varepsilon)}\cap [0, r]\}.
	\]
	Then $\eta_0>1$ and for all $n\geq 1$,
	\[
	x_j(n+m) =T_j(x(n+m-1))=x_j(n+m-1)f_j(x(n+m-1)) \geq \eta_0 x_j(n+m-1).
	\]
	This leads to $x_j(n+m) \geq \eta^n_0 x_j(m)\to +\infty$ as $n\to +\infty$, a contradiction to the boundedness of $\{x(n)\}$.
	
	In the latter case, from Step 1 we see that $\lim_{n\to +\infty}y_k(n) =0$. Without loss of generality, we may assume that $0<y_k<q_k$, where $Q_k= q_ke_k$. Since the whole trajectory $\gamma(y)$ is contained in $\omega(x)$ and $\omega(x)$ is compact, from Step 1 we derive that $\omega(y) \subset \omega(x)\cap(\cap_{\ell=1}^{k-1}\pi_{\ell})\cap\Sigma$. Let
	\[
	\eta_1 = \max\{f_k(u): u\in [0, r]\}
	\]
	and take a small $\varepsilon\in (0, \frac{y_k}{3\eta_1})$. Since the set 
	\[
	S = \{z\in (\cap_{\ell=1}^{k-1}\pi_{\ell})\cap\Sigma: \varepsilon \leq z_k\leq y_k\} 
	\]
	is compact and $\omega(z)\subset (\cap_{\ell=1}^k\pi_{\ell})\cap\Sigma$ for all $z\in S$ from Step 1, by continuous dependence there is a $\delta \in (0, \min\{\frac{\varepsilon}{2}, \frac{1}{6}y_k\})$ such that 
	\begin{eqnarray}
	\forall u\in O(S, \delta)\cap [0, r], &\exists z\in S,& \exists n_1(z)\in\N, \; \textup{such that}\; z_k(n_1) < \frac{\varepsilon}{2}, \nonumber \\
	& & \forall n\in\{0, \ldots, n_1\}, \|u(n)-z(n)\|< \frac{\varepsilon}{2}. \label{e17}
	\end{eqnarray}
	For this $\delta$, there is an $m_0\in\N$ such that $x(n)\in O((\cap_{\ell=1}^{k-1}\pi_{\ell})\cap\Sigma, \delta)$ for all $n\geq m_0$. As $\emptyset \not= \omega(y)\subset\omega(x)\cap(\cap_{\ell=1}^k\pi_{\ell})\cap \Sigma$, there is an $m_1 \geq m_0$ such that $0 <x_k(m_1)< \varepsilon$. Then
	\[
	x_k(m_1 +1) = T_k(x(m_1)) = x_k(m_1)f_k(x(m_1)) \leq \eta_1 \varepsilon < \frac{1}{3}y_k.
	\]
	Thus, either (i) $x_k(m_1+1) < \varepsilon$ or (ii) $x(m_1+1)\in O(S, \delta)$. In case (i), by taking $m_2 = m_1 +1$ we have $x_k(m_2+1) < \frac{1}{3}y_k$. In case (ii), by (\ref{e17}) there is a $z\in S$ and $n_1\in\N$ such that $z_k(n_1)< \frac{1}{2}\varepsilon$ and $\|x(m_1+1 +n) -z(n)\|<\frac{1}{2}\varepsilon$ for all $0 \leq n\leq n_1$. By the positive invariance of $(\cap_{\ell=1}^{k-1}\pi_{\ell})\cap\Sigma$, for any $z\in (\cap_{\ell=1}^{k-1}\pi_{\ell})\cap\Sigma$, we have $z(n)\in (\cap_{\ell=1}^{k-1}\pi_{\ell})\cap\Sigma$ for all $n\in\N$. As (\ref{e13}) implies that $(\cap_{\ell=1}^{k-1}\pi_{\ell})\cap\Sigma$ is above $\Gamma_k$ and $S\subset (\cap_{\ell=1}^{k-1}\pi_{\ell})\cap\Sigma$, $z_k(n)$ is monotone nonincreasing in $n$ for each $z\in S$. Then,
	\[
	x_k(m_1+1+n) < z_k(n) + \frac{1}{2}\varepsilon \leq z_k + \frac{1}{2}\varepsilon < x_k(m_1+1) +\varepsilon <\frac{2}{3}y_k, 0\leq n\leq n_1
	\]
	and $x_k(m_1+1+n_1) < z_k(n_1) + \frac{1}{2}\varepsilon < \varepsilon$. Take $m_2 = m_1+1+n_1$. In either (i) or (ii), we see that $x_k(n) < \frac{2}{3}y_k$ for all $m_1 \leq n\leq m_2$. Repeating the above process we obtain $x_k(n) < \frac{2}{3}y_k$ for all $n\geq m_1$, a contradiction to $y\in\omega(x)$.
	
	The contradictions in both cases above show that $\omega(x) \subset (\cap_{\ell=1}^k\pi_{\ell})\cap\Sigma$. Then (\ref{e15}) follows.  
	
	Finally, if $k = N-1$, we have $\lim_{n\to +\infty}x(n) = Q_N$ for all $x\in C$ with $x_N>0$. As $Q_N$ is above $\Gamma_j$ for all $j\in \{1, \ldots, N-1\}$ and $\frac{\partial f_N}{\partial x_N}(Q_N) <0$, every eigenvalue of the Jacobian matrix $DT(Q_N)$ is in the interval $(0, 1)$. Thus, $Q_N$ is globally asymptotically stable. 
\end{proof}

Note that the statement of Theorem \ref{The3.2} used the natural ascending order of numbers for the species. Obviously, the statement is still true after a permutation from ascending order of numbers. 

\begin{corollary}\label{Cor3.1}
	Assume that (\ref{e6'}) and the conditions of Theorem \ref{The2.2} hold. Assume also the existence of a permutation $p: I_N \to I_N$ and an integer $k\in I_N\setminus\{N\}$ such that $\frac{\partial f_{p(i)}}{\partial x_{p(i)}}(Q_{p(i)})<0$ for all $i\in \{1, \ldots, k\}$ and
	\begin{equation}\label{e18}
	\forall j\in\{i+1, \ldots, N\}, (\cap^{i-1}_{\ell=1}\pi_{p(\ell)})\cap \Gamma_{p(i)}\cap[0, r]\; \textit{is strictly below}\; \Gamma_{p(j)}.
	\end{equation}
	Then the $p(i)$th species is dominated for all $i\in \{1, \ldots, k\}$. In addition, if $k = N-1$ and $\frac{\partial f_{p(N)}}{\partial x_{p(N)}}(Q_{p(N)})<0$ then the $p(N)$th species is dominant and the $p(N)$th axial fixed point $Q_{p(N)}$ is globally asymptotically stable in $C$.
\end{corollary}

\section{Some Examples}\label{Sec4}
In this section, we apply our results obtained in sections 2 and 3 to some known models as examples. All these models fit well into our system (\ref{e1}) for maps $T: C\to C$ of the form (\ref{e2}), where the sign of each entry $\frac{\partial f_i(x)}{\partial x_j}$ of the Jacobian $Df$ is completely determined by the corresponding entry $a_{ij}$ of a constant matrix $A = (a_{ij})_{N\times N}$:
\[
\forall i, j\in I_N, \forall x\in C, \frac{\partial f_i(x)}{\partial x_j} = -\sigma_{ij}(x)a_{ij}, \sigma_{ij}(x)>0.
\]
All entries of the matrix $A$ are assumed positive in most of the references cited here either due to their particular meaning in the original model or due to convenience of theoretical analysis by using available results such as Theorem \ref{The1.1}, Theorem \ref{The1.2} and their variations mentioned in section 1. Under such an assumption, each system models the population dynamics of a community of $N$ competing species where the population of the $j$th species directly affects the growth rate of the population of the $i$th species in a negative way as $a_{ij}>0$.  

Unfortunately, as far as the author knows, not much has been found about such systems modelling competing species where the population of the $j$th species affects the growth rate of the population of the $i$th species in a negative way, directly or indirectly, due to $a_{ij}>0$ or $a_{ij}=0$. With the help of our Theorem \ref{The2.1} and Theorem \ref{The2.2}, we are now able to deal with these models under the relaxed assumption:
\[
\forall i, j\in I_N, a_{ij} \geq 0, a_{ii} >0.
\] 
Since our results obtained in section 3 and this section below are all based on the assumptions of Theorem \ref{The2.2}, if $\frac{\partial f_i(x)}{\partial x_j}=0$ ($a_{ij}=0$ for the models below) for at least one pair of indices $i, j$ at some point $x\in[0, r]$, then these results are not achievable by using Theorem \ref{The1.2} and its variations as their conditions are not fully met. This demonstrates that the class of systems to which Theorem \ref{The2.2} is applicable is broader than that for Theorem \ref{The1.2} and its variations. Hence, our main results are a significant improvement of those available in literature. 

\subsection{The competitive Leslie-Gower models}
The competitive Leslie-Gower models are system (\ref{e1}) for maps $T: C\to C$ of the form (\ref{e2}): $T_i(x) = x_if_i(x)$, where
\begin{equation}\label{e19}
\forall i, j\in I_N, f_i(x) = \frac{c_i}{1+\sum^N_{k=1}a_{ik}x_k}, c_i>1, a_{ij}\geq 0, a_{ii} >0.
\end{equation}
Under the condition that $a_{ij}>0$ for all $i, j\in I_N$, Jiang and Niu \cite{JiNi2} have shown that each Leslie-Gower model admits a carrying simplex. 

I. Following the same lines as those in \cite{JiNi2}, we check that each Leslie-Gower model with $a_{ij}\geq 0$ and $a_{ii}>0$ has a unique modified carrying simplex $\Sigma$ by Theorem \ref{The2.2}. Indeed, for $x\in X_i$, $f_i(x) = 1$ if and only if $x_i= \frac{c_i-1}{a_{ii}} = q_i$, so $T$ restricted to $X_i$ has a unique fixed point $Q_i= q_ie_i$. For $i, j\in I_N$,
\[
\frac{\partial f_i}{\partial x_j} = -\frac{c_ia_{ij}}{[1+\sum^N_{k=1}a_{ik}x_k]^2},
\]
so $\frac{\partial f_i}{\partial x_j} \leq 0$ and $\frac{\partial f_i}{\partial x_i} <0$. Also, for all $x\in C$,
\[
\forall i\in I_N, f_i(x) +\sum_{j=1}^Nx_j \frac{\partial f_i}{\partial x_j} = \frac{c_i}{[1+\sum^N_{k=1}a_{ik}x_k]^2} >0.
\]
By Remark 2.2 (c), conditions (i)--(iii) of Theorem \ref{The2.2} and (\ref{e6'}) are all met for any $r\gg q$.

The surfaces $\Gamma_i$ are now $(N-1)$-dimensional planes in $C$:
\[
\forall i\in I_N, \Gamma_i = \{x\in C: a_{i1}x_1+ \cdots + a_{iN}x_N = c_i-1\}.
\]

II. If for some $i\in I_N$, the following inequalities hold:
\begin{equation}\label{e20}
\forall j\in I_N\setminus\{i\}, \forall k\in I_N, a_{jk}(c_i-1)<a_{ik}(c_j-1), 
\end{equation}
then the intersection point of $\Gamma_i$ with each positive half axis $X_k$ is below $\Gamma_j$ for every $j\in I_N\setminus\{i\}$. So $\Gamma_i$ is strictly below $\Gamma_j$ in $C$ for all $j\in I_N\setminus\{i\}$. By Theorem \ref{The3.1} (a), the $i$th species is dominated. 

III. If for some $i\in I_N$, 
\begin{equation}\label{e21}
\forall j\in I_N\setminus\{i\}, \forall k\in I_N, \; \textup{either}\; a_{ik}=a_{jk} = 0 \; \textup{or} \; a_{ik}(c_j-1)<a_{jk}(c_i-1), 
\end{equation}
then either $X_k$ is parallel to both $\Gamma_i$ and $\Gamma_j$ or the intersection point of $\Gamma_j$ with $X_k$ is below $\Gamma_i$, so $\Gamma_i$ is strictly above $\Gamma_j$ for all $j\in I_N\setminus\{i\}$. By Theorem \ref{The3.1} (b), the $i$th species is dominant and the fixed point $Q_i$ is globally asymptotically stable.

IV. Note that (\ref{e21}) is a sufficient condition for $\Gamma_i$ to be strictly above $\Gamma_j$ for all $j\in I_N\setminus\{i\}$ in $C$. But the condition in Theorem \ref{The3.1} (b) only requires the relationship of such planes restricted to $[0, r]$. So (\ref{e21}) is much stronger than the requirement of Theorem \ref{The3.1} (b). For example, let us consider the three-dimensional Leslie-Gower model with
\begin{equation}\label{e22}
f_1(x) =\frac{2}{1+ x_1 + 0.25x_3}, f_2(x)=\frac{2}{1+2x_1+x_2+0.2x_3}, f_3(x) = \frac{2}{1+2x_1+ x_3}.
\end{equation}
Clearly, $e_1$, $e_2$ and $e_3$ are the axial fixed points. Take $r = (1.1, 1.1, 1.1)\gg (1, 1, 1) =q$. The intersection points of $\Gamma_3$ with the $X_1$ and $X_3$ are $(0.5, 0, 0)$ and $e_3$ respectively. As $f_1(0.5, 0, 0) = \frac{2}{1.5} >1$ and $f_1(e_3) = \frac{2}{1.25}>1$, both $(0.5, 0, 0)$ and $e_3$ are below $\Gamma_1$. Since $X_2$ is parallel to both $\Gamma_1$ and $\Gamma_3$, $\Gamma_1$ is strictly above $\Gamma_3$ in $\R^3_+$. The intersection points of $\Gamma_2$ with the axes are $(0.5, 0, 0)$, $e_2$ and $(0, 0, 5)$. We know that $(0.5, 0, 0)$ is below $\Gamma_1$ already. As $f_1(e_2) = 2>1$ but $f_1(0, 0, 5) = \frac{2}{2.25}<1$, $e_2$ is below $\Gamma_1$ but $(0, 0, 5)$ is above $\Gamma_1$. So $\Gamma_1$ is not above $\Gamma_2$ on $\R^3_+$ and (\ref{e21}) is not met. However, restricted to $[0, r]$, $\Gamma_2$ intersects one of the edges of $[0, r]$ at $(0.39, 0, 1.1)$ and $f_1(0.39, 0, 1.1) = \frac{2}{1.665}>1$. So $\Gamma_2\cap\pi_2\cap [0, r]$ is strictly below $\Gamma_1$. This, together with $e_2$ below $\Gamma_1$, implies that $\Gamma_1\cap [0, r]$ is strictly above $\Gamma_2$. By Theorem \ref{The3.1} (b), the first species is dominant and the fixed point $e_1$ is globally asymptotically stable.

V. Now suppose the following inequalities hold:
\begin{eqnarray}
\forall i\in I_N\setminus\{N\}, \forall j, k &\in& \{i+1, \ldots, N\}, \label{e23} \\ 
a_{ji}(c_i-1) &<& a_{ii}(c_j-1),\; a_{(i+1)k}(c_i-1) <a_{ik}(c_{i+1}-1). \nonumber
\end{eqnarray}
Then, for each $i\in I_N\setminus\{N\}$, the intersection point of $\Gamma_i$ with $X_i$ is below $\Gamma_j$ for all $j\in\{i+1, \ldots, N\}$ and $\Gamma_i\cap(X_{i+1}\times\cdots\times X_N)$ is strictly below $\Gamma_{i+1}$. Thus, $\Gamma_{N-1}\cap(X_{N-1}\times X_N)$ is strictly below $\Gamma_N$, $\Gamma_{N-2}\cap (X_{N-2}\times X_{N-1}\times X_N)$ is strictly below $\Gamma_{N-1}$ and $\Gamma_N$, $\ldots$, $\Gamma_i\cap(X_i\times\cdots\times X_N)$ is strictly below $\Gamma_j$ for all $j\in\{i+1, \ldots, N\}$. By Theorem \ref{The3.2}, the $N$th species is dominant and $Q_N = \frac{c_N-1}{a_{NN}}e_N$ is globally asymptotically stable.

\subsection{The generalised competitive Atkinson-Allen models}
The generalised competitive Atkinson-Allen models are systems (\ref{e1}) for maps $T: C\to C$ of the form (\ref{e2}): $T_i(x) = x_if_i(x)$, where
\begin{equation}\label{e24}
\forall i, j\in I_N, f_i(x) = c_i+ \frac{(1+u_i)(1-c_i)}{1+\sum^N_{k=1}a_{ik}x_k}, 0<c_i <1, u_i>0, a_{ij}\geq 0, a_{ii} >0.
\end{equation}
Under the condition that $a_{ij}>0$ for all $i, j\in I_N$, Gyllenberg et al \cite{GyJiNiYa2} have shown that each such model admits a carrying simplex. 

I. Following the same lines as those in \cite{GyJiNiYa2}, we check that each generalised Atkinson-Allen model with $a_{ij}\geq 0$ and $a_{ii}>0$ has a unique modified carrying simplex $\Sigma$ by Theorem \ref{The2.2}. Indeed, for $x\in X_i$, $f_i(x) = 1$ if and only if $x_i= \frac{u_i}{a_{ii}} = q_i$, so $T$ restricted to $X_i$ has a unique fixed point $Q_i= q_ie_i$. For $i, j\in I_N$,
\[
\frac{\partial f_i}{\partial x_j} = -\frac{(1+u_i)(1-c_i)a_{ij}}{[1+\sum^N_{k=1}a_{ik}x_k]^2},
\]
so $\frac{\partial f_i}{\partial x_j} \leq 0$ and $\frac{\partial f_i}{\partial x_i} <0$. Also, for all $x\in C$,
\[
\forall i\in I_N, f_i(x) +\sum_{j=1}^Nx_j \frac{\partial f_i}{\partial x_j} = c_i+\frac{(1+u_i)(1-c_i)}{[1+\sum^N_{k=1}a_{ik}x_k]^2} >0.
\]
By Remark 2.2 (c), conditions (i)--(iii) of Theorem \ref{The2.2} and (\ref{e6'}) are all met for any $r\gg q$.

The surfaces $\Gamma_i$ are now $(N-1)$-dimensional planes in $C$:
\begin{equation}\label{e24'}
\forall i\in I_N, \Gamma_i = \{x\in C: a_{i1}x_1+ \cdots + a_{iN}x_N = u_i\}.
\end{equation}

II. If for some $i\in I_N$, the following inequalities hold:
\begin{equation}\label{e25}
\forall j\in I_N\setminus\{i\}, \forall k\in I_N, a_{jk}u_i<a_{ik}u_j, 
\end{equation}
then the intersection point of $\Gamma_i$ with each positive half axis $X_k$ is below $\Gamma_j$ for every $j\in I_N\setminus\{i\}$. So $\Gamma_i$ is strictly below $\Gamma_j$ in $C$ for all $j\in I_N\setminus\{i\}$. By Theorem \ref{The3.1} (a), the $i$th species is dominated. 

III. If for some $i\in I_N$, 
\begin{equation}\label{e26}
\forall j\in I_N\setminus\{i\}, \forall k\in I_N, \; \textup{either}\; a_{ik}=a_{jk} = 0 \; \textup{or} \; a_{ik}u_j <a_{jk}u_i, 
\end{equation}
then either $X_k$ is parallel to both $\Gamma_i$ and $\Gamma_j$ or the intersection point of $\Gamma_j$ with $X_k$ is below $\Gamma_i$, so $\Gamma_i$ is strictly above $\Gamma_j$ for all $j\in I_N\setminus\{i\}$. By Theorem \ref{The3.1} (b), the $i$th species is dominant and the fixed point $Q_i$ is globally asymptotically stable.

IV. Note that (\ref{e26}) is much stronger than the requirement of Theorem \ref{The3.1} (b). Similar to (\ref{e22}), we can easily construct a three-dimensional generalised Atkinson-Allen model as an example which fails (\ref{e26}) but satisfies the condition of Theorem \ref{The3.1} (b).

V. Now suppose the following inequalities hold:
\begin{equation}\label{e27}
\forall i\in I_N\setminus\{N\}, \forall j, k\in\{i+1, \ldots, N\},\; a_{ji}u_i < a_{ii}u_j,\; a_{(i+1)k}u_i <a_{ik}u_{i+1}.
\end{equation}
Then, for each $i\in I_N\setminus\{N\}$, the intersection point of $\Gamma_i$ with $X_i$ is below $\Gamma_j$ for all $j\in\{i+1, \ldots, N\}$ and $\Gamma_i\cap(X_{i+1}\times\cdots\times X_N)$ is strictly below $\Gamma_{i+1}$. Thus, $\Gamma_{N-1}\cap(X_{N-1}\times X_N)$ is strictly below $\Gamma_N$, $\Gamma_{N-2}\cap (X_{N-2}\times X_{N-1}\times X_N)$ is strictly below $\Gamma_{N-1}$ and $\Gamma_N$, $\ldots$, $\Gamma_i\cap(X_i\times\cdots\times X_N)$ is strictly below $\Gamma_j$ for all $j\in\{i+1, \ldots, N\}$. By Theorem \ref{The3.2}, the $N$th species is dominant and $Q_N = \frac{u_N}{a_{NN}}e_N$ is globally asymptotically stable.

VI. The standard Atkinson-Allen models are systems (\ref{e1}) for maps $T: C\to C$ of the form (\ref{e2}): $T_i(x) = x_if_i(x)$, where
\begin{equation}\label{e28}
\forall i, j\in I_N, f_i(x) = c+ \frac{2(1-c)}{1+\sum^N_{k=1}a_{ik}x_k}, 0<c <1, a_{ij}\geq 0, a_{ii} >0.
\end{equation}
Note that $f$ defined by (\ref{e28}) is a special case of (\ref{e22}) with $c_i = c$ and $u_i = 1$ for all $i\in I_N$. Thus, the results obtained above for generalised Atkinson-Allen models can be applied to the standard Atkinson-Allen models with simplified conditions ($u_i, u_j$ replaced by 1 in (\ref{e25})--(\ref{e27})). For these models with $a_{ij}>0$ for all $i, j\in I_N$, see \cite{DiWaYa}, \cite{JiNi1} and the references therein for further results.

\subsection{The competitive Ricker models}
The competitive Ricker models are systems (\ref{e1}) for maps $T: C\to C$ of the form (\ref{e2}): $T_i(x) = x_if_i(x)$, where
\begin{equation}\label{e29}
\forall i, j\in I_N, f_i(x) = \exp{\left[u_i\left(1-\sum^N_{k=1}a_{ik}x_k\right)\right]}, u_i>0, a_{ij}\geq 0, a_{ii} >0.
\end{equation}
Under the conditions that $a_{ij}>0$ for all $i, j\in I_N$ and
\begin{equation}\label{e30}
\forall i\in I_N, u_i< a_{ii}/\sum^N_{j=1}a_{ij}, \; \textup{or} \; \forall i\in I_N, u_i < 1/\sum^N_{j=1}\frac{a_{ij}}{a_{jj}},
\end{equation}
Gyllenberg et al \cite{GyJiNiYa1} have shown that each such model admits a carrying simplex. 

I. We check that, under (\ref{e30}), each Ricker model with $a_{ij}\geq 0$ and $a_{ii}>0$ has a unique modified carrying simplex $\Sigma$ by Theorem \ref{The2.2}. (i) For $x\in X_i$, $f_i(x) = 1$ if and only if $x_i= \frac{1}{a_{ii}} = q_i$, so $T$ restricted to $X_i$ has a unique fixed point $Q_i= q_ie_i$. (ii) For $i, j\in I_N$,
\[
\frac{\partial f_i}{\partial x_j} = -u_ia_{ij}f_i(x),
\]
so $\frac{\partial f_i}{\partial x_j} \leq 0$ and $\frac{\partial f_i}{\partial x_i} <0$. (iii) For all $x\in [0, q]$, we have
\[
\forall i\in I_N, f_i(x) +\sum_{j=1}^Nx_j \frac{\partial f_i}{\partial x_j} = f_i(x)[1 - u_i\sum^N_{j=1}a_{ij}x_j] \geq f_i(x)[1 - u_i\sum^N_{j=1}\frac{a_{ij}}{a_{jj}}] >0,
\]
or
\[
\forall i\in I_N, f_i(x) +x_i\sum_{j=1}^N \frac{\partial f_i}{\partial x_j} = f_i(x)[1 - u_ix_i\sum^N_{j=1}a_{ij}] \geq f_i(x)[1 - \frac{u_i}{a_{ii}}\sum^N_{j=1}a_{ij}] >0.
\]
By Remark 2.2 (c), conditions (i)--(iii) of Theorem \ref{The2.2} and (\ref{e6'}) are all met for any $r\gg q$. Then, by Theorem \ref{The2.2}, each Ricker model with (\ref{e30}) has a modified carrying simplex $\Sigma$.

The surfaces $\Gamma_i$ are now $(N-1)$-dimensional planes in $C$:
\[
\forall i\in I_N, \Gamma_i = \{x\in C: a_{i1}x_1+ \cdots + a_{iN}x_N = 1\}.
\]

II. If for some $i\in I_N$, the following inequalities hold:
\begin{equation}\label{e31}
\forall j\in I_N\setminus\{i\}, \forall k\in I_N, a_{jk}<a_{ik}, 
\end{equation}
then the intersection point of $\Gamma_i$ with each positive half axis $X_k$ is below $\Gamma_j$ for every $j\in I_N\setminus\{i\}$. So $\Gamma_i$ is strictly below $\Gamma_j$ in $C$ for all $j\in I_N\setminus\{i\}$. By Theorem \ref{The3.1} (a), the $i$th species is dominated. 

III. If for some $i\in I_N$, 
\begin{equation}\label{e32}
\forall j\in I_N\setminus\{i\}, \forall k\in I_N, \; \textup{either}\; a_{ik}=a_{jk} = 0 \; \textup{or} \; a_{ik} <a_{jk}, 
\end{equation}
then either $X_k$ is parallel to both $\Gamma_i$ and $\Gamma_j$ or the intersection point of $\Gamma_j$ with $X_k$ is below $\Gamma_i$, so $\Gamma_i$ is strictly above $\Gamma_j$ for all $j\in I_N\setminus\{i\}$. By Theorem \ref{The3.1} (b), the $i$th species is dominant and the fixed point $Q_i$ is globally asymptotically stable.

IV. Note that (\ref{e32}) is much stronger than the requirement of Theorem \ref{The3.1} (b). Similar to (\ref{e22}), we can easily construct a three-dimensional Ricker model as an example which fails (\ref{e32}) but satisfies the condition of Theorem \ref{The3.1} (b).

V. Now suppose the following inequalities hold:
\begin{equation}\label{e33}
\forall i\in I_N\setminus\{N\}, \forall j, k\in\{i+1, \ldots, N\},\; a_{ji} < a_{ii},\; a_{(i+1)k} <a_{ik}.
\end{equation}
Then, for each $i\in I_N\setminus\{N\}$, the intersection point of $\Gamma_i$ with $X_i$ is below $\Gamma_j$ for all $j\in\{i+1, \ldots, N\}$ and $\Gamma_i\cap(X_{i+1}\times\cdots\times X_N)$ is strictly below $\Gamma_{i+1}$. Thus, $\Gamma_{N-1}\cap(X_{N-1}\times X_N)$ is strictly below $\Gamma_N$, $\Gamma_{N-2}\cap (X_{N-2}\times X_{N-1}\times X_N)$ is strictly below $\Gamma_{N-1}$ and $\Gamma_N$, $\ldots$, $\Gamma_i\cap(X_i\times\cdots\times X_N)$ is strictly below $\Gamma_j$ for all $j\in\{i+1, \ldots, N\}$. By Theorem \ref{The3.2}, the $N$th species is dominant and $Q_N = \frac{1}{a_{NN}}e_N$ is globally asymptotically stable.

\subsection{General competitive models with plane nullclines}
In \cite{Hou}, the competitive models given by system (\ref{e1}) for maps $T: C\to C$ of the form (\ref{e2}) $T_i(x) = x_if_i(x)$ are considered, where
\begin{equation}\label{e34}
\forall i\in I_N, f_i(x) = G_i((Ax)_i)
\end{equation}
with
\begin{equation}\label{e35}
A=\left(\begin{array}{llll} a_{ii} & a_{12} & \cdots & a_{1N} \\ a_{21} & a_{22} & \cdots & a_{2N} \\ \cdots & \cdots & \cdots &\cdots \\ a_{N1} &a_{N2} & \cdots & a_{NN} \end{array} \right) 
\end{equation}
satisfying $a_{ii}>0$ and $a_{ij}\geq 0$, $(Ax)_i$ denoting the $i$th component of $Ax$. Assume that the functions $G_i\in C^1(\R_+, \R_+)$ satisfy the following conditions:
\begin{itemize}
	\item[(a1)] Each $G_i$ is positive and strictly decreasing with $G_i(u_i) = 1$ and $G_i'(u_i) <0$ for some $u_i>0$.
	\item[(a2)] For $x\in C$ and each $i\in I_N$, $\frac{\partial T_i}{\partial x_i} >0$ for $0\leq x_i \leq \frac{u_i}{a_{ii}}=q_i$.
\end{itemize}
Then each nullcline surface $\Gamma_i$ is a hyperplane given by (\ref{e24'}). Under (a1), (a2) and another condition, criteria are established in \cite{Hou} for global stability of a fixed point by geometric method of using the relative positions of the nullcline planes in $[0, q]$. The matrix $M(x)$ defined by (\ref{e4}) is $M(x) = -\textup{diag}(\frac{x_iG_i'((Ax)_i)}{G_i((Ax)_i)})A$.

I. Assume that
\begin{equation}\label{e36}
\forall x\in [0, q], \rho(M(x))<1.
\end{equation}
We check that each such model has a unique modified carrying simplex $\Sigma$ by Theorem \ref{The2.2}. (i) For $x\in X_i$, $G_i((Ax)_i) = 1$ if and only if $(Ax)_i = u_i$, i.e. $x_i= \frac{u_i}{a_{ii}} = q_i$, so $T$ restricted to $X_i$ has a unique fixed point $Q_i= q_ie_i$. (ii) For $i, j\in I_N$,
\[
\frac{\partial G_i((Ax)_i)}{\partial x_j} = a_{ij}G_i'((Ax)_i).
\]
By $a_{ij}\geq 0$, $a_{ii}>0$ and (a1), $\frac{\partial G_i((Ax)_i)}{\partial x_j} \leq 0$ and $G_i((Ax)_i)$ is strictly decreasing in $x_i$. Condition (iii) of Theorem \ref{The2.2} follows from (\ref{e36}). Then, by Theorem \ref{The2.2}, each model with (\ref{e34}) has a modified carrying simplex $\Sigma$. 

By the same reasoning as that given for generalised Atkinson-Allen models, we obtain the following conclusions.

II. If (\ref{e25}) holds, then the $i$th species is dominated.

III. If (\ref{e26}) holds, then the $i$th species is dominant and $Q_i$ is globally asymptotically stable.

IV. If (\ref{e27}) holds, then the $N$th species is dominant and $Q_N$ is asymptotically stable.

\section{Proof of the main theorems}\label{Sec5}
In this section, we aim at providing complete proofs for Theorem \ref{The2.1}, Theorem \ref{The2.2} and Corollary \ref{Cor2.2}. Although some of the ideas used here are credited to \cite{Her} and \cite{WaJi1, WaJi2}, for readers' convenience we present an independent proof rather than citing some lemmas and theorems and modifying their proofs bit by bit. However, this does not mean that the proofs are trivial modifications from those in the references. Actually, the author's main contribution in this paper is the sharp observation that the system permits a modifies carrying simplex if the retrotone property for $T$ is relaxed to weakly retrotone, which leads to the dramatic relaxation of the conditions of Theorem \ref{The1.2} to those of Theorem \ref{The2.2}. To prove these results, in addition to inheriting some techniques shown in Lemma \ref{Lem5.2}, Lemma \ref{Lem5.3}, the main part of the proof of Theorem \ref{The2.1} and a small part of the proof of Theorem \ref{The2.2}, the author's own methods and techniques are reflected in Lemma \ref{Lem5.1}, Lemma \ref{Lem5.4}, Lemma \ref{Lem5.5}, the main part of the proof of Theorem \ref{The2.2} and Corollary \ref{Cor2.2}. 

\begin{lemma}\label{Lem5.1}
	Assume that $T$ satisfies the conditions of Theorem \ref{The2.1}. Then, for any $x\in [0, r]$, $[0, T(x)]\subset T([0, r])$.
\end{lemma}
\begin{proof}
	By Remark 2.1 (b) we know that $T$ is a homeomorphism from $[0, r]$ to $T([0, r])$. Thus, $T$ maps an open set of $[0, r]$ to an open set of $T([0, r])$. Clearly, the set
	\[
	[0, r) = \{x\in [0, r]: 0\leq x\ll r\}
	\] 
	is open in $[0, r]$, so $T([0, r))$ is also open in $T([0, r])$. We first show that
	\begin{equation}\label{5e1}
	\forall x\in [0, r), [0, T(x)] \subset T([0, r)).
	\end{equation}
	Suppose (\ref{5e1}) is not true. Then, for some $x\in [0, r)\setminus\{0\}$, there is a $y$ satisfying $0<y<T(x)$ but $y\not\in T([0, r))$. Since $T([0, r))$ is open and $T(x)\in T([0, r))$, there is an $s_0\in [0, 1)$ such that $y(s) = y +s(T(x)-y) \in T([0, r))$ for $s\in (s_0, 1]$ but $y(s_0) = y+s_0(T(x)-y)\not\in T([0, r))$. As $T$ is weakly retrotone, we have $0< z(s) = T^{-1}(y(s)) < x$ for $s\in (s_0, 1)$ and $z(s_1)<z(s_2)$ for any $s_0 <s_1<s_2\leq 1$. Thus, $\lim_{s\to s_0+}z(s)$ exists. Define $z(s_0) = \lim_{s\to s_0+}z(s)$. Then $z(s_0) \in [0, x]\subset [0, r)$ so $T(z(s_0)) \in T([0, r))$. By continuity of $T$,
	\[
	T(z(s_0)) = \lim_{s\to s_0+}T(z(s)) = \lim_{s\to s_0+}y(s) = y(s_0),
	\]
	a contradiction to $y(s_0)\not\in T([0, r))$. This contradiction shows the truth of (\ref{5e1}).
	
	Now we show that $[0, T(x)] \subset T([0, r])$ for all $x\in [0, r]$. This is true by (\ref{5e1}) if $x\in [0, r)$, so we suppose $x\in [0, r]\setminus [0, r)$. Then $T(x) \in T([0, r]\setminus [0, r))$ and $x(s) = sx \in [0, r)$ for all $s\in [0, 1)$ with $\lim_{s\to 1-}x(s) = x$. Moreover, by (\ref{5e1}), $[0, T(x(s))] \subset T([0, r)) \subset T([0, r])$ for all $s\in [0, 1)$. For each $y\in [0, T(x)]$, if $y\in [0, T(x(s))]$ for some $s\in [0, 1)$ then $y\in T([0, r])$; if $y\not\in [0, T(x(s))]$ for any $s\in [0, 1)$ then there is an increasing sequence $\{s_n\}\subset [0, 1)$ with $s_n\uparrow 1$ and a sequence $\{y_n\}$ with $y_n\in [0, T(x(s_n))]\subset T([0, r))$ such that $\lim_{n\to \infty}y_n =y$. So $y\in \overline{T([0, r))}$. But since $T([0, r))\subset T([0, r])$, $T([0, r))$ is open and $T([0, r])$ is closed, we have $\overline{T([0, r))} \subset T([0, r])$ so $y\in T([0, r])$. This shows that $[0, T(x)]\subset T([0, r])$ for all $x\in [0, r]$.
\end{proof}

For any $x\in [0, r]$, we denote the image of $x$ under $(T^{-1})^k$ by $x(-k)$ if $(T^{-1})^k(x)=T^{-k}(x)$ exists.

\begin{lemma}\label{Lem5.2}
	Assume that the conditions of Theorem \ref{The2.1} hold. Suppose $x\in [0, r]\setminus\{0\}$ such that $(T^{-1})^k(x)$ exists and $x(-k)\in [0, r]$ for all $k\in\N$. Then, for any $y\in [0, r]$ with $y<x$ and $x-y\in \dot{C}_I$ for some nonempty $I\subset I_N$, $y(-k)$ exists in $[0, r]$ for all $k\in \N$ and
	\begin{equation}\label{5e2}
	\forall i\in I, \lim_{k\to \infty}y_i(-k) = 0.
	\end{equation}
\end{lemma}
\begin{proof}
	By the existence of $x(-1) \in [0, r]$ we have $x = T(T^{-1}(x)) = T(x(-1)) \in T([0, r])$. Thus, by Lemma \ref{Lem5.1}, $[0, x] \subset T([0, r])$. As $y\in [0, r]$ and $y<x$, we have $y\in T([0, r])$ so $y(-1)$ exists and $y(-1)\in [0, r]$. It then follows from the weak retrotone property of $T$ that $y(-1) <x(-1)$ and $y_i(-1) < x_i(-1)$ for all $i\in I$. If $y(-k) = (T^{-1})^k(y)$ exists, $y(-k) <x(-k)$ and $y_i(-k) < x_i(-k)$ for all $i\in I$ and some $k\in\N$, by the same reasoning as above we obtain the existence of $y(-k-1) = (T^{-1})^{k+1}(y)$, $y(-k-1) <x(-k-1)$ and $y_i(-k-1) < x_i(-k-1)$ for all $i\in I$. By induction, we see the existence of $y(-k)\in [0, r]$ with $y(-k)<x(-k)$ and $y_i(-k)<x_i(-k)$ for all $i\in I$ and all $k\in\N$.
	
	To prove (\ref{5e2}) by contradiction, we suppose the existence of $i\in I$ such that $0< y_i<x_i$ and $y_i(-k)\not\to 0$ as $k\to\infty$. As $x(-k), y(-k)\in [0, r]$ for all $k\in\N$ and $[0, r]$ is compact, we can select a subsequence $\{\sigma(k)\}\subset \{k\}$ such that 
	\[
	\lim_{k\to\infty}x(-\sigma(k)) = \bar{x}, \lim_{k\to\infty}y(-\sigma(k)) =\bar{y}, \bar{y}_i >0.
	\]
	By $0<y_i<x_i$ we have $0<y_i(-k)<x_i(-k)$ for all $k\in\N$. Now define
	\[
	\Delta(k) = \frac{y_i(-k)}{x_i(-k)}, k\in\N.
	\]
	Then
	\[
	1> \Delta(k) = \frac{T_i(y(-k-1))}{T_i(x(-k-1))} = \frac{y_i(-k-1)f_i(y(-k-1))}{x_i(-k-1)f_i(x(-k-1))} = \Delta(k+1)\frac{f_i(y(-k-1))}{f_i(x(-k-1))}.
	\]
	By condition (iii) of Theorem \ref{The2.1}, $\frac{f_i(y(-k-1))}{f_i(x(-k-1))}>1$. So
	\[
	\forall k\in\N, 0< \Delta(k+1) <\Delta(k) <1.
	\]
	This shows the existence of a $\beta\in [0, 1)$ such that $\lim_{k\to\infty}\Delta(k) = \beta$. In particular,
	\[
	1> \beta = \lim_{k\to\infty}\Delta(\sigma(k)) =\lim_{k\to\infty}\frac{y_i(-\sigma(k))}{x_i(-\sigma(k))} = \frac{\bar{y}_i}{\bar{x}_i} >0,
	\]
	so $0 <\bar{y}_i =\beta\bar{x}_i < \bar{x}_i$. By continuity of $T$, $T(x(-\sigma(k)))\to T(\bar{x})$ and $T(y(-\sigma(k)))\to T(\bar{y})$ as $k\to\infty$. Thus,
	\[
	\beta = \lim_{k\to\infty}\Delta(\sigma(k)-1) =\lim_{k\to\infty}\frac{y_i(-\sigma(k)+1)}{x_i(-\sigma(k)+1)} = \lim_{k\to\infty}\frac{T_i(y(-\sigma(k)))}{T_i(x(-\sigma(k)))} = \frac{T_i(\bar{y})}{T_i(\bar{x})}. 
	\]
	From this we obtain $T_i(\bar{y}) = \beta T_i(\bar{x}) < T_i(\bar{x})$. As $y(-k) <x(-k)$ for all $k\in\N$, we have $T(y(-\sigma(k))) = y(-\sigma(k)+1) < x(-\sigma(k)+1) = T(x(-\sigma(k)))$ and, as $k\to\infty$, $T(\bar{y}) < T(\bar{x})$. By condition (iii) of Theorem \ref{The2.1}, we have $\frac{f_i(\bar{y})}{f_i(\bar{x})} >1$. As $\lim_{k\to\infty}\frac{f_i(y(-\sigma(k)))}{f_i(x(-\sigma(k)))} = \frac{f_i(\bar{y})}{f_i(\bar{x})}$ and $\frac{f_i(y(-k))}{f_i(x(-k))} >1$ for all $k\in\N$, there is an $\eta >1$ such that $\frac{f_i(y(-\sigma(k)))}{f_i(x(-\sigma(k)))} \geq \eta$ for all $k\in\N$. Thus,
	\[
	\Delta(\sigma(k)-1) =\frac{T_i(y(-\sigma(k)))}{T_i(x(-\sigma(k)))} =\Delta(\sigma(k))\frac{f_i(y(-\sigma(k)))}{f_i(x(-\sigma(k)))} \geq \eta\Delta(\sigma(k)).
	\]
	From this and $\sigma(k)-1 \geq \sigma(k-1)$ for $k>1$ we obtain
	\[
	\forall k>1, \Delta(\sigma(k)) \leq \frac{1}{\eta}\Delta(\sigma(k)-1) \leq \frac{1}{\eta}\Delta(\sigma(k-1)).
	\]
	This implies that
	\[
	\forall k\in\N, \beta < \Delta(\sigma(k+1)) \leq \frac{1}{\eta^k}\Delta(\sigma(1)).
	\]
	Letting $k\to\infty$, we obtain $\beta =0$, a contradiction to $\beta = \frac{\bar{y}_i}{\bar{x}_i} >0$. This contradiction shows the truth of (\ref{5e2}).
\end{proof}

\begin{lemma}\label{Lem5.3}
	Suppose the existence of $x, y\in [0, r]\setminus\{0\}$ with support $I(x) = I(y) \subset I_N$ satisfying $x(k) \leq y(k)$ for all $k\in\N$. Then, under the conditions of Theorem \ref{The2.1}, 
	\begin{equation}\label{5e3}
	\lim_{k\to\infty}(y(k)-x(k)) = 0.
	\end{equation}
\end{lemma}
\begin{proof}
	If $x(k_0) = y(k_0)$ for some $k_0\in\N$ then $y(k) = x(k)$ for all $k\geq k_0$ so (\ref{5e3}) holds. Now assume that $x(k)<y(k)$ for all $k\in\N$. For each $i\in I(x)$, if there is a $k_1\in\N$ such that $x_i(k_1) = y_i(k_1)$, then we must have $x_i(k_1+1) = y_i(k_1+1)$, for the inequality $T_i(x(k_1))= x_i(k_1+1) < y_i(k_1+1) = T_i(y(k_1))$ and condition (ii) of Theorem \ref{The2.1} would imply $x_i(k_1)<y_i(k_1)$. Thus, $x_i(k) = y_i(k)$ for all $k\geq k_1$ so $\lim_{k\to\infty}(y_i(k) - x_i(k)) =0$. Now suppose for a fixed $i\in I(x)$,
	\[
	\forall k\in\N, 0<x_i(k)<y_i(k).
	\]
	Define $\delta(k) = \frac{x_i(k)}{y_i(k)}$ for all $k\in\N$. Then
	\[
	1 > \delta(k+1) = \frac{T_i(x(k))}{T_i(y(k))}=\delta(k) \frac{f_i(x(k))}{f_i(y(k))}.
	\]
	As $ \frac{f_i(x(k))}{f_i(y(k))}>1$ by condition (iii) of Theorem \ref{The2.1}, $\{\delta(k)\}$ is a positive increasing sequence bounded above by 1. If $\lim_{k\to\infty}\delta(k) = 1$ then
	\[
	y_i(k)-x_i(k) = y_i(k)[1-\delta(k)] \to 0 \; (k\to\infty).
	\]
	Suppose $\lim_{k\to\infty}\delta(k) = \beta$ for some $\beta\in (0, 1)$. If $\lim_{k\to\infty}(y_i(k)-x_i(k)) \not= 0$, there must be a subsequence $\{\sigma(k)\}\subset \{k\}$ such that
	\[
	\lim_{k\to\infty}x(\sigma(k)) = \bar{x}, \lim_{k\to\infty}y(\sigma(k)) = \bar{y}, \bar{x} < \bar{y}, \bar{x}_i<\bar{y}_i.
	\]
	Then,
	\begin{equation}\label{5e4}
	1> \delta(\sigma(k+1))\geq \delta(\sigma(k)+1) = \frac{T_i(x(\sigma(k)))}{T_i(y(\sigma(k)))} = \delta(\sigma(k)) \frac{f_i(x(\sigma(k)))}{f_i(y(\sigma(k)))}.
	\end{equation}
	By condition (iii) of Theorem \ref{The2.1} again, $\frac{f_i(x(\sigma(k)))}{f_i(y(\sigma(k)))}>1$. As
	\[
	T(x(\sigma(k))) = x(\sigma(k)+1) < y(\sigma(k)+1) = T(y(\sigma(k))),
	\]
	$\lim_{k\to\infty}T(x(\sigma(k))) = T(\bar{x})$ and $\lim_{k\to\infty}T(y(\sigma(k))) = T(\bar{y})$, we have $T(\bar{x}) \leq T(\bar{y})$. If $T_i(\bar{x}) < T_i(\bar{y})$ then $T(\bar{x}) < T(\bar{y})$. By condition (iii) of Theorem \ref{The2.1}, we obtain $\frac{f_i(\bar{x})}{f_i(\bar{y})} >1$. If $T_i(\bar{x}) = T_i(\bar{y})$, then
	\[
	1= \frac{T_i(\bar{x})}{T_i(\bar{y})} = \frac{\bar{x}_i}{\bar{y}_i}\frac{f_i(\bar{x})}{f_i(\bar{y})} , \; \frac{f_i(\bar{x})}{f_i(\bar{y})} =\frac{\bar{y}_i}{\bar{x}_i} >1.
	\]
	Therefore, there is an $\eta>1$ such that $\frac{f_i(x(\sigma(k)))}{f_i(y(\sigma(k)))}\geq \eta$ for all $k\in\N$. Then, from (\ref{5e4}) we obtain
	\[
	\delta(\sigma(k+1)) \geq \eta \delta(\sigma(k)) \geq \eta^k\delta(\sigma(1)) \to +\infty \; (k\to\infty),
	\]
	a contradiction to $\delta(k) <1$. This contradiction shows the conclusion (\ref{5e3}).
\end{proof}

Under the assumptions of Theorem \ref{The2.1}, $T([0, r]) \subset [0, r]$. By Remark 2.1 (b), $T: [0, r] \to T([0, r])$ is a homeomorphism, so $T$ maps open sets to open sets and closed sets to closed sets. As 
\[
\forall n\in\N, T^{n+1}([0, r]) \subset T^n([0, r])
\]
and $[0, r]$ is compact, $T^n([0, r])$ is compact for all $n\in\N$. From Remark 2.1 (e) we know that 0 is a repellor with basin of repulsion $\mathcal{B}(0)\subset [0, r] $.

\begin{lemma}\label{Lem5.4}
	Assume that the conditions of Theorem \ref{The2.1} hold. Let
	\begin{equation}\label{5e5}
	A_0 = \cap^{\infty}_{n=0} T^n([0, r]).
	\end{equation}
	Then $A_0$ is nonempty, compact, invariant and $A_0 = \overline{\mathcal{B}(0)}$.
\end{lemma}
\begin{proof}
	That $A_0\not= \emptyset$ is obvious as $0, q_ie_i$ and all fixed points of $T$ are in $A_0$. As each $T^n([0, r])$ is compact and any nonempty intersection of compact sets is compact, by (\ref{5e5}) $A_0$ is compact. The invariance of $A_0$ follows from (\ref{5e5}) and $T([0, r]) \subset [0, r]$. Clearly, by (\ref{5e5}) we see that $A_0$ is the largest invariant set of $T$ in $[0, r]$. As $\mathcal{B}(0)$ is an open subset of $[0, r]$ and invariant, we have $\mathcal{B}(0)\subset A_0$. To show that $A_0 =\overline{\mathcal{B}(0)}$, we take an arbitrary point $x\in A_0 \setminus\mathcal{B}(0)$ and show that $x\in \overline{\mathcal{B}(0)}$. This is trivial if $x=0$ as $0\in\overline{\mathcal{B}(0)}$. If $x\not= 0$ then there is a nonempty $I\subset I_N$ such that $x\in \dot{C}_I$. Moreover, $u_s = sx\ll_I x$ for all $s\in (0, 1)$. By the invariance of $A_0$, $x(-k)$ exists in $A_0$ for all $k\in\N$. Then, by Lemma \ref{Lem5.2}, $u_s(-k)$ exists in $[0, r]$ for all $s\in (0, 1)$ and all $k\in\N$ and
	$\lim_{k\to\infty}u_s(-k) =0$. Thus, $u_s \in \mathcal{B}(0)$ for all $s\in (0, 1)$. Since $x = \lim_{s\to 1-}u_s$, we have $x\in \overline{\mathcal{B}(0)}$. 
\end{proof}

With the help of Lemmas \ref{Lem5.1}--\ref{Lem5.4} we are now in a position to prove Theorem \ref{The2.1}.

\begin{proof}[Proof of Theorem \ref{The2.1}]
	Let $\Sigma = \overline{\mathcal{B}(0)} \setminus(\{0\}\cup \mathcal{B}(0))$. We verify that $\Sigma$ is a modified carrying simplex. Clearly, $\Sigma \subset [0, r]\setminus\{0\}$ and $\Sigma\not=\emptyset$ as all the nontrivial fixed points are in $\Sigma$. From Lemma \ref{Lem5.4} we see that $\overline{\mathcal{B}(0)}$ is compact and invariant. As $\{0\}\cup \mathcal{B}(0)$ is open and invariant, $\Sigma$ is compact and invariant. That $T: \Sigma\to\Sigma$ is a homeomorphism follows from $T$ being a homeomorphism from $[0, r]$ to $T([0, r])$.
	
	To show that $\Sigma$ is homeomorphic to $\Delta^{N-1}$ by radial projection, we define a map $m: \Delta^{N-1}\to\Sigma$ as follows. For each $x\in \Delta^{N-1}$, as $\lambda x\in\mathcal{B}(0)$ for sufficiently small $\lambda >0$ and $\mathcal{B}(0)$ is open, there is a unique $\lambda_0 = \lambda_0(x)>0$ such that $\lambda x\in\mathcal{B}(0)$ for all $0<\lambda <\lambda_0$ but $\lambda_0x\not\in\mathcal{B}(0)$. Since $\lim_{\lambda\to\lambda_0}\lambda x = \lambda_0x \not= 0$, we have $\lambda_0x\in\overline{\mathcal{B}(0)}$ so $\lambda_0x\in\Sigma$. We claim that $\lambda x\not\in A_0$ for $\lambda >\lambda_0$, where $A_0$ is given by (\ref{5e5}). Indeed, if there is a $\lambda_1 > \lambda_0$ such that $u =\lambda_1 x\in A_0$, then $u(-k) \in A_0$ for all $k\in\N$. By Lemma \ref{Lem5.2}, we would have $\lambda_0 x \in \mathcal{B}(0)$, a contradiction to $\lambda_0x\not\in\mathcal{B}(0)$. Thus,
	\[
	\forall x\in \Delta^{N-1}, \Sigma\cap\{\lambda x: \lambda>0\} = \{\lambda_0(x)x\}.
	\]
	Then the map $m: \Delta^{N-1} \to\Sigma$ defined by $m(x) =\lambda_0(x)x$ is a bijection. The map $m$ is a homeomorphism if $m$ and $m^{-1}$ are continuous.
	
	To show that $m$ is continuous, we need only show that $\lambda_0: \Delta^{N-1}\to\R_+$ is continuous. Suppose $\lambda_0$ is not continuous at a point $x_0\in\Delta^{N-1}$, i.e. $\lim_{x\to x_0}\lambda_0(x) \not= \lambda_0(x_0)$. Since $\lambda_0$ is obviously bounded, there is a sequence $\{x_k\} \subset \Delta^{N-1}$ such that
	\[
	x_k\to x_0 \;\textup{and} \; \lambda_0(x_k) \to \mu \not= \lambda_0(x_0)\; \textup{as}\; k\to\infty. 
	\]
	Then $\{m(x_k)\}\subset \Sigma$ and $m(x_k) = \lambda_0(x_k)x_k \to \mu x_0$ as $k\to\infty$. Since $\Sigma$ is compact, we have $\mu x_0\in\Sigma$. This contradicts $\Sigma\cap\{\lambda x: \lambda >0\} = \{\lambda_0(x_0)x_0\} \not= \{\mu x_0\}$. This contradiction shows the continuity of $m$ on $\Delta^{N-1}$.
	
	To show that $m^{-1}: \Sigma\to\Delta^{N-1}$ is continuous, since the continuity of $m$ implies that $\Sigma$ is a continuous surface, for each $y\in\Sigma$, there is a unique $\mu=\mu(y)>0$ such that $\mu(y)y\in \Delta^{N-1}$ so that $m^{-1}(y) = \mu(y)y$. Then the continuity of $m^{-1}$ follows from showing the continuity of $\mu: \Sigma\to\R_+$ by the same technique as above. Therefore, $\Sigma$ is homeomorphic to $\Delta^{N-1}$ by radial projection.
	
	Next, we show that for each $x\in [0, r]\setminus\{0\}$, if $x$ is above $\Sigma$ then $\omega(x)\subset \Sigma$; if $x$ is below $\Sigma$ then there is a $y\in\Sigma$ with support $I(y) = I(x)$ such that 
	\begin{equation}\label{5e6}
	\lim_{k\to +\infty}[x(k)-y(k)] = 0.  
	\end{equation} 
	
	Now suppose $x$ is above $\Sigma$. By Lemma \ref{Lem5.4} we have $\omega(x) \subset A_0$. As $\mathcal{B}(0)$ does not contain any positive limit point and $0\not\in\omega(x)$, we must have $\omega(x)\subset \Sigma$.
	
	Next, suppose $x$ is below $\Sigma$ with support $I(x)\subset I_N$. As $\overline{\mathcal{B}(0)}=\{0\}\cup\mathcal{B}(0)\cup\Sigma$, we have $\Sigma^- =\{0\}\cup\mathcal{B}(0)$. By $x>0$, we must have $x\in\mathcal{B}(0)$. Define sets
	\[
	\forall k\in\N, U(k, x) = \{y\in\Sigma: I(y) = I(x), x(k) <y(k)\}.
	\]  
	Note that $T(\mathcal{B}(0)) = \mathcal{B}(0)$ and $T(\Sigma) = \Sigma$. For each fixed $k\in\N$, $x(k)$ is below $\Sigma$ so there is $\mu_0>1$ such that $\mu_0 x(k)\in\Sigma$ but $\mu x(k) \in\mathcal{B}(0)$ for $1\leq \mu <\mu_0$. Taking $y= T^{-k}(\mu_0 x(k))$ we have $y\in\Sigma$, $I(y) = I(x)$ and $y(k) = \mu_0 x(k) >x(k)$. Thus, $y\in U(k, x)$ so $U(k, x)\not=\emptyset$. For each $z\in U(k+1, x)$ we have
	\[
	T(x(k)) = x(k+1) < z(k+1) = T(z(k)).
	\]
	As $T$ is weakly retrotone, we must have $x(k) < z(k)$ so $z\in U(k, x)$. This shows that
	\[
	\forall k\in\N, U(k+1, x) \subset U(k, x).
	\]
	From the definition we see that each $U(k, x)$ is compact. Then $\emptyset\not= \cap^{\infty}_{k=0} U(k, x) \subset \Sigma$. Taking any $y\in \cap^{\infty}_{k=0} U(k, x)$ we obtain $x(k)<y(k)$ for all $k\in\N$ so (\ref{5e6}) follows from Lemma \ref{Lem5.3}.
	
	So far we have proved that $\Sigma = \overline{\mathcal{B}(0)}\setminus(\{0\}\cup\mathcal{B}(0))$ is a modified carrying simplex. Now for each $p\in\Sigma$ and every $q\in[0, r]\setminus\{0\}$ with $q<p$, by Lemma \ref{Lem5.2} we know that $q(-k)\in [0, r]$ exists for all $k\in\N$ and $\lim_{k\to +\infty}q_i(-k) = 0$, so $\alpha(q) \subset \pi_i$, for any $i\in I_N$ with $q_i<p_i$.
	
	Finally, we show the uniqueness of the modified carrying simplex $\Sigma$. Suppose we have another modified carrying simplex $\Sigma_1 \not= \Sigma$. Then, on a half line starting from the origin we have two distinct points $p\in\Sigma$ and $q\in\Sigma_1$ so there is a positive number $\lambda \not= 1$ such that $p = \lambda q$. Clearly $p$ and $q$ have the same support $I\subset I_N$ so we have either $p\ll_I q$ or $q\ll_I p$. In the first case, by Lemma \ref{Lem5.2} we would have $\alpha(p) =\{0\}$, a contradiction to $\alpha(p) \subset \Sigma$ as $0\not\in \Sigma$. In the second case, by Lemma \ref{Lem5.2} again we would have $\alpha(q) = \{0\}$, a contradiction to $\alpha(q) \subset\Sigma_1$ as $0\not\in\Sigma_1$. This shows that $\Sigma$ is the unique modified carrying simplex.
\end{proof}

To prove Theorem \ref{The2.2}, we need the following lemma.

\begin{lemma}\label{Lem5.5}
	Let $U$ be a small neighbourhood of $[0, r]$ and $T\in C^1(U, U)$. Assume that the Jacobian matrix $DT(x)$ is invertible on $[0, r]$ with $(DT(x))^{-1} = (t_{ij})$. If
	\[
	\forall x\in U, \forall i, j\in I_N, \; t_{ii}(x) > 0\; \textup{and}\; t_{ij}(x) \geq 0,
	\]
	then $T$ from $U$ to $T(U)$ is one-to-one and is weakly retrotone on $[0, r]$.
\end{lemma}
\begin{proof}
	Since $DT(x)$ is continuous on $U$ and invertible on $[0, r]$, there is a small neighbourhood $U_1\subset U$ of $[0, r]$ such that $DT(x)$ is invertible on $U_1$. Without loss of generality, we assume that $U_1=U$. By the inverse function theorem, $T$ from $U$ to $T(U)$ is one-to-one and invertible. Moreover, $T^{-1}$ on $T(U)$ is differentiable. As $g(u)= u = T^{-1}(T(u))$ for $u\in U$, by the chain rule of differentiation we have
	\[
	I = Dg(u) = [D(T^{-1})(T(u))][DT(u)], 
	\]  
	so 
	\begin{equation}\label{5e7}
	D(T^{-1})(T(u)) = (DT(u))^{-1} =(t_{ij}(u)). 
	\end{equation}
	Now for any $x, y\in [0, r]$ with $T(x)<T(y)$ and $T(y)-T(x) \in \dot{C}_I$ for some nonempty $I\subset I_N$, we have
	\begin{eqnarray*}
		y-x &=& T^{-1}(T(y)) - T^{-1}(T(x)) \\
		&=& T^{-1}(T(x)+s(T(y)-T(x)))|^1_0 \\
		&=& \int^1_0 \frac{d}{ds}T^{-1}(T(x)+s(T(y)-T(x))) ds \\
		&=& \int^1_0D(T^{-1})(T(x)+s(T(y)-T(x))) ds (T(y)-T(x)).
	\end{eqnarray*}
	By the assumption on the entries of $(DT(u))^{-1}$ and (\ref{5e7}), the the diagonal entries of the matrix $D(T^{-1})(T(x)+s(T(y)-T(x)))$ are positive and other entries are nonnegative. As the matrix 
	\[
	\int^1_0D(T^{-1})(T(x)+s(T(y)-T(x))) ds 
	\]
	maintains the same feature as $D(T^{-1})(T(x)+s(T(y)-T(x)))$, for each $i\in I$, $T(y)-T(x)>0$ and $T_i(y)-T_i(x) >0$ imply $y-x >0$ and $y_i-x_i>0$. Thus, $T$ on $[0, r]$ is weakly retrotone.
\end{proof}

\begin{proof}[Proof of Theorem \ref{The2.2}]
	We need only show that conditions (ii) and (iii) of Theorem \ref{The2.2} imply conditions of (ii) and (iii) of Theorem \ref{The2.1}. Since
	\[
	DT(x) = \textup{diag}(f_1(x), \ldots, f_N(x))(I - M(x)),
	\]
	where $I$ is the identity matrix and $M(x)$ is given by (\ref{e4}), if $\rho(M(x))<1$ by (iii) of Theorem \ref{The2.2}, then $DT(x)$ is invertible with
	\begin{eqnarray*}
		(DT(x))^{-1} &=& (I - M(x))^{-1}\textup{diag}\bigl(\frac{1}{f_1(x)}, \ldots, \frac{1}{f_N(x)}\bigr) \\
		&=& \bigl(I+\sum^{\infty}_{k=1} M^k(x)\bigr)\textup{diag}\bigl(\frac{1}{f_1(x)}, \ldots, \frac{1}{f_N(x)}\bigr).
	\end{eqnarray*}
	From this it is clear that each diagonal entry of $(DT(x))^{-1}$ is positive and other entries are nonnegative. Then condition (ii) of Theorem \ref{The2.1} follows from Lemma \ref{Lem5.5}. 
	
	Now suppose $\rho(\tilde{M}(x))<1$ holds instead of $\rho(M(x))<1$, where $\tilde{M}(x)$ is given by (\ref{E4}). If $x\gg 0$, then
	\[
	\textup{diag}(\frac{1}{x_1}, \ldots, \frac{1}{x_N})DT(x)\textup{diag}(x_1, \ldots, x_N) = \textup{diag}(f_1(x), \ldots, f_N(x))(I - \tilde{M}(x)),
	\]
	so $DT(x)$ is invertible with
	\begin{eqnarray*}
		(DT(x))^{-1} &=& \textup{diag}(x_1, \ldots, x_N)(I - \tilde{M}(x))^{-1}\textup{diag}\bigl(\frac{1}{x_1f_1(x)}, \ldots, \frac{1}{x_Nf_N(x)}\bigr) \\
		&=& \textup{diag}(x_1, \ldots, x_N)\bigl(I+\sum^{\infty}_{k=1} \tilde{M}^k(x)\bigr)\textup{diag}\bigl(\frac{1}{x_1f_1(x)}, \ldots, \frac{1}{x_Nf_N(x)}\bigr).
	\end{eqnarray*}
	From this we see that each diagonal entry of $(DT(x))^{-1}$ is positive and other entries are nonnegative. Then condition (ii) of Theorem \ref{The2.1} follows from Lemma \ref{Lem5.5}. 
	
	If $\rho(\tilde{M}(x))<1$ holds with $x\not\gg 0$, then there is a proper subset $J_1\subset I_N$ as the support of $x$. Without loss of generality, we assume that $J_1 = \{1, \ldots, k\}$ for some positive integer $k<N$ (as we can always rearrange the order of the components). Let $J_2 = \{k+1, \ldots, N\}$ and $U = \textup{diag}(x_1, \ldots, x_k, 1, \ldots, 1)$. Then
	\[
	U^{-1}DT(x)U = \textup{diag}(f_1(x), \ldots, f_N(x))\left(\begin{array}{cc} I_1-M_1(x) & -M_2(x) \\ 0 & I_2 \end{array}\right),
	\]
	where $I_1$ and $I_2$ are $k\times k$ and $(N-k) \times(N-k)$ identity matrices respectively, and
	\begin{eqnarray*}
		M_1(x) &=& \left(-\frac{x_j}{f_i(x)}\frac{\partial f_i}{\partial x_j}(x)\right)_{k\times k}, \;\textup{for}\; i, j\in J_1, \\
		M_2(x) &=& \left(-\frac{1}{f_i(x)}\frac{\partial f_i}{\partial x_j}(x)\right)_{k\times (N-k)}, \;\textup{for}\; i\in J_1, j\in J_2.
	\end{eqnarray*}
	Note that
	\[
	\tilde{M}(x) = \left(\begin{array}{cc} M_1(x) & 0 \\ M_3(x) & 0 \end{array}\right),
	\]
	where $M_3(x)$ is an $(N-k)\times k$ matrix. Then $\rho(M_1(x)) =\rho(\tilde{M}(x)) <1$, so $DT(x)$ is invertible with
	\begin{eqnarray*}
		(DT(x))^{-1} &=& U\left(\begin{array}{cc} (I - M_1(x))^{-1} & (I - M_1(x))^{-1}M_2(x)\\ 0 & I_2 \end{array}\right) \\
		& & \times\textup{diag}\bigl(\frac{1}{f_1(x)}, \ldots, \frac{1}{f_N(x)}\bigr)U^{-1}. 
	\end{eqnarray*}
	As $(I - M_1(x))^{-1}= I +\sum^{\infty}_{n=1} M_1^n(x)$ with positive diagonal entries and nonnegative other entries, each diagonal entry of $(DT(x))^{-1}$ is positive and other entries are nonnegative. Then condition (ii) of Theorem \ref{The2.1} follows from Lemma \ref{Lem5.5}. 
	
	For any $x, y\in [0, r]$, if $T(x) < T(y)$ and $T(y)-T(x) \in \dot{C}_J$ for some $J\subset I_N$, by the weakly retrotone property of $T$ we have $x<y$ and $x_j<y_j$ for all $j\in J$. By condition (ii) of Theorem \ref{The2.2}, each $f_i$ is nonincreasing in every $x_j$ but strictly decreasing in $x_i$ for $x\in [0, r]$. Then we have $f(x)>f(y)$ and $f_j(x)>f_j(y)$ for all $j\in J$, so condition (iii) of Theorem \ref{The2.1} holds.
	
	Finally, we check that $[0, r]$ is positively invariant. Note that Remark 2.1 (b) and (c) do not reply on the positive invariance of $[0, r]$. Then, for each $x\in [0, r]$ and every $i\in I_N$, by (ii) of Theorem \ref{The2.2} we have $T_i(x) \leq T_i(x_ie_i)$. By Remark 2.1 (c), $T_i(x_ie_i)$ is increasing for $x_i\in [0, r_i]$, so
	\[
	T_i(x) \leq T_i(x_ie_i) \leq T_i(r_ie_i) = r_if_i(r_ie_i) < r_i.
	\]
	This shows that $T(x) \ll r$ and $T([0, r]) \subset [0, r]$.
\end{proof}

\begin{proof}[Proof of Corollary \ref{Cor2.2}]
	For any bounded set $B\subset [0, r]\setminus\{0\}$ with $\overline{B}\subset [0, r]\setminus\{0\}$, there is a small $\delta >0$ such that $O(0, \delta)\cap \overline{B} =\emptyset$ and $O(0, \delta)\cap [0, r]$ is strictly below $\Sigma$. Since $0$ is a repellor, $\overline{\mathcal{B}(0)}$ is invariant by Lemma \ref{Lem5.4}, $[0, r]$ is positively invariant, and $T$ from $[0, r]$ to $T([0, r])$ is a homeomorphism by Remark 2.1 (b), for $\delta >0$ small enough the set $[0, r]\setminus O(0, \delta)$ is positively invariant with $\overline{B} \subset [0, r]\setminus O(0, \delta)$ and $O(0, \delta)\cap [0, r]$ is strictly below $\Sigma$. Then, for each $n\in \N$, $T^n([0, r]\setminus O(0, \delta))$ is compact and 
	\[
	\Sigma \subset T^{n+1}([0, r]\setminus O(0, \delta)) \subset T^n([0, r]\setminus O(0, \delta)).
	\]
	From this follows
	\[
	\Sigma \subset  \bigcap_{n=0}^{\infty} T^n([0, r]\setminus O(0, \delta)).
	\]
	We claim that
	\begin{equation}\label{e101}
	\Sigma =  \bigcap_{n=0}^{\infty} T^n([0, r]\setminus O(0, \delta)).
	\end{equation}
	Indeed, from Lemma \ref{Lem5.4} we know that 
	\[
	\cap_{n=0}^{\infty} T^n([0, r]\setminus O(0, \delta)) \subset A_0 = \overline{\mathcal{B}(0)} = \{0\}\cup \mathcal{B}(0)\cup\Sigma.
	\]
	If (\ref{e101}) is not true, then there is a point $p\in (\cap_{n=0}^{\infty} T^n([0, r]\setminus O(0, \delta)))\setminus \Sigma$, so $T^{-n}(p)\in [0, r]\setminus O(0, \delta)$ for all $n\in\N$. This shows that $\lim_{n\to\infty}T^{-n}(p) \not= 0$. On the other hand, however, as $p\not\in \{0\}\cup\Sigma$, we must have $p\in \mathcal{B}(0)$ so $\lim_{n\to\infty}T^{-n}(p) = 0$. This contradiction shows the truth of (\ref{e101}).
	
	Now from (\ref{e101}) we see that $\Sigma$ attracts the points of $[0, r]\setminus O(0, \delta)$ uniformly. As $B\subset [0, r]\setminus O(0, \delta)$, $\Sigma$ attracts the points of $B$ uniformly. Therefore, $\Sigma$ is a global attractor in $[0, r]\setminus\{0\}$ under the conditions of Theorem \ref{The2.1} or Theorem \ref{The2.2}.
	
	Under the additional condition (\ref{e6'}), for any bounded set $B\subset C\setminus\{0\}$ with $\overline{B}\subset C\setminus\{0\}$, from Remark 2.3 we know the existence of an integer $k>0$ such that $T^k(\overline{B})\subset [0, r]$. By the definition of $T$, $T^k(x)=0$ if and only if $x=0$ on $C$. As $0\not\in\overline{B}$, $0\not\in T^k(\overline{B})$ so $T^k(\overline{B})\subset [0, r]\setminus\{0\}$. From the previous paragraph we know that $\Sigma$ attracts the points of $T^k(\overline{B})$ uniformly. Thus, $\Sigma$ attracts the points of $B$ uniformly. Hence,  $\Sigma$ is a global attractor in $C\setminus\{0\}$. 	
\end{proof}

\section{Conclusion}\label{Sec6}
We have so far considered the discrete dynamical system (\ref{e1}) with the maps $T$ defined by (\ref{e2}). Recall that the current available carrying simplex theory is about the existence of an $(N-1)$-dimensional surface that is a compact invariant set attracting all the points in $C\setminus\{0\}$. With the existing concept of a carrying simplex for the system and the available criteria on existence of carrying simplex as the main concern of this paper, we have successfully achieved our goal of extending this theory to a broader class of systems: We have first defined the concept of a modified carrying simplex, which is a slight relaxation from the concept of a carrying simplex and is still an $(N-1)$-dimensional surface that is a compact invariant global attractor of the system in $C\setminus\{0\}$. We then have established our criteria for existence and uniqueness of a modified carrying simplex.

In comparison with the existing criteria for existence of a carrying simplex, our criteria for existence and uniqueness of a modified carrying simplex have the following main virtue: Instead of requiring all the entries of the Jacobian $Df(x)$ to be negative for all $x\in [0, q]$, we only require each entry of $Df(x)$ to be nonpositive and each $f_i(x)$ to be strictly decreasing in $x_i$. Thus, we have significantly reduced the cost of having an $(N-1)$-dimensional surface as a compact invariant global attractor of the system in $C\setminus\{0\}$. In other words, our criteria can be applied to a broader class of systems as competitive models.  

The significance of the carrying simplex theory lies in that the global dynamics of the system in $C$ can be described by the dynamics on the modified carrying simplex $\Sigma$. As one application of this theory, we have investigated vanishing species and dominance of one species over others. Assuming the existence of a modified carrying simplex, we have obtained sufficient geometric conditions for one or more species to die out. We have also obtained conditions for one species to dominate all others and one axial fixed point to be globally asymptotically stable.

Above all, with our theorems for modified carrying simplex $\Sigma$, we have laid the foundation for exploring the global dynamics of the system. We expect future research work will be flourishing based on modified carrying simplex.

\vskip 3 mm
\noindent \textbf{Open Problem} Suppose system (\ref{e1}) with $T$ defined by (\ref{e2}) satisfies the conditions of Theorem \ref{The2.2}, so the system permits a modified carrying simplex $\Sigma$. Is it possible to construct a sequence $\{T^{[k]}\}$ satisfying the following conditions?
\begin{itemize}
	\item[(i)] For each integer $k>0$, the map $T^{[k]}$ from $[0, r]$ to $C$ has the form (\ref{e2}).
	\item[(ii)] Each $T^{[k]}$ on $[0, r]$ meets the requirements of Theorem \ref{The1.2}, so system (\ref{e1}) with $T$ replaced by $T^{[k]}$ permits a carrying simplex $\Sigma^{[k]}$.
	\item[(iii)] As $k\to\infty$, $T^{[k]}(x)\to T(x)$ uniformly for $x\in [0, r]$
	\item[(iv)] As $k\to\infty$, $\Sigma^{[k]} \to\Sigma$ in the following sense:
	\[
	\forall \varepsilon >0, \exists K>0, \forall k\geq K, \Sigma^{[k]} \subset O(\Sigma, \varepsilon).
	\]
\end{itemize}
If the answer is YES, then our Theorem \ref{The2.2} can be viewed as the result of a limit process from Theorem \ref{The1.2}, i.e. system (\ref{e1}) with (\ref{e2}) satisfying the conditions of Theorem \ref{The2.2} can be approximated by systems satisfying the conditions of Theorem \ref{The1.2}.

\section*{Acknowledgements} The author consulted Professor Stephen Baigent on this topic and is grateful for his encouragement of writing up this paper.
The author is also grateful to the referees and editors for their comments and suggestions adopted in this version of the paper.


\medskip
\medskip

\end{document}